\documentclass[11pt]{article}

\usepackage{latexsym,amssymb}

\makeatletter \@addtoreset{equation}{section}

\newtheorem{thm}{Theorem}[section]
\newtheorem{prop}[thm]{Proposition}
\newtheorem{lem}[thm]{Lemma}
\newtheorem{cor}[thm]{Corollary}
\newtheorem{defin}[thm]{Definition}

\def\enu#1{\newline\makebox[5mm][l]{\rm(#1)}}
\def\benu#1#2{\newline\makebox[#1 mm][l]{\rm(#2)}}

\def\bp{\noindent{\it Proof.}\ }
\def\bpp#1{\noindent{\it Proof of #1.}\ }
\def\ep{\nopagebreak\newline\mbox{\ }\hfill\rule{2mm}{2mm}}
\def\epp{\nopagebreak\mbox{\ }\hfill\rule{2mm}{2mm}}

\def\ad{{\rm ad}\,}
\renewcommand{\Im}{{\rm Im}}
\def\mod{{\rm mod}\,}
\def\span{{\rm span}}
\def\supp{{\rm supp}\,}
\def\Tr{{\rm Tr}}

\def\7#1{{\mathbb #1}}

\def\A{{\cal A}}
\def\B{{\cal B}}
\def\C{{\cal C}}
\def\F{{\cal F}}
\def\I{{\cal I}}

\def\iso{{\raisebox{1mm}{$\sim$}\hspace{-4.2mm}\to}}

\def\2{{\frac{1}{2}}}
\def\ii{{\frac{i}{2}}}

\def\D{{\Delta}}
\def\sl2{{U_q({\mathfrak su}_2)}}

\begin{document}

\title{\bf The Martin Boundary of a Discrete Quantum Group}

\author{Sergey Neshveyev \& Lars Tuset}
\date{}

\maketitle

\begin{abstract}
We consider the Markov operator $P_\phi$ on a discrete quantum
group given by convolution with a $q$-tracial state $\phi$. In the
study of harmonic elements $x$, $P_\phi(x)=x$, we define the
Martin boundary $A_\phi$. It is a separable C$^*$-algebra carrying
canonical actions of the quantum group and its dual. We establish
a representation theorem to the effect that positive harmonic
elements correspond to positive linear functionals on $A_\phi$.
The C$^*$-algebra $A_\phi$ has a natural time evolution, and the
unit can always be represented by a KMS state. Any such state
gives rise to a u.c.p. map from the von Neumann closure of
$A_\phi$ in its GNS representation to the von Neumann algebra of
bounded harmonic elements, which is an analogue of the Poisson
integral. Under additional assumptions this map is an isomorphism
which respects the actions of the quantum group and its dual. Next
we apply these results to identify the Martin boundary of the dual
of $SU_q(2)$ with the quantum homogeneous sphere of Podle{\'s}.
This result extends and unifies previous results by Ph.~Biane and
M.~Izumi.
\end{abstract}

\section*{Introduction}

Hopf algebras play an important role in quantum probability as
they provide a natural setting for generalizations of L\'{e}vy
processes \cite{SCH}. The convolution operators on Hopf algebras
also lead to one of the simplest non-trivial examples of
non-commutative Markov processes. In a series of papers
\cite{B1,B2,B3,B4} Biane studied such operators on the dual $\hat
G$ of a compact group $G$. He considered the group von Neumann
algebra $L(G)$ of $G$ with comultiplication
$\D(\lambda_g)=\lambda_g\otimes\lambda_g$ and looked at the
convolution operator $P_\phi =(\phi\otimes\iota)\D$, where $\phi$
is a normal tracial state on $L(G)$. One is particularly
interested in $P_\phi$-harmonic elements, that is, the elements
$x$ affiliated with $L(G)$ such that $P_\phi(x)=x$. A significant
part of the theory of random walks on abelian groups \cite{R,S}
can be generalized to this context. For example, one can prove
that the constants are the only bounded harmonic
elements~\cite{B1}. This result also follows from a generalization
of the Choquet-Deny theorem, which states that extremal normalized
harmonic elements are exponentials, that is, positive group-like
elements \cite{B2}. In the classical theory the set of extremal
elements constitutes a part of the boundary of an appropriate
compactification of $\hat G$. This boundary, coined the Martin
boundary, is constructed by completing $\hat G$ with respect to a
metric depending on the asymptotic behavior of a function $K$
called the Martin kernel. Having described the set of extremal
harmonic elements, one could ask for a geometric realization of
this set as a boundary of $\hat G$. The boundary should then be
understood in the sense of non-commutative geometry, so it should
be a unital C$^*$-subalgebra of $L(G)/C^* (G)$. This problem was
solved for $SU(2)$ in \cite{B4}, where Biane introduced a
non-commutative analogue of the Martin kernel, proved that the
corresponding boundary of $\widehat{SU(2)}$ is the $2$-dimensional
sphere, and showed that this sphere could be naturally identified
with the set of harmonic exponentials. Thus he obtained an
analogue of the Ney-Spitzer theorem which describes the Martin
compactification of $\7Z^n$, see e.g. \cite{W}. Note, however,
that for this result Biane assumes $\|\phi\|<1$, so he considers
sub-Markov operators for which all harmonic elements are
unbounded. If $\phi$ is a state, then there are, in fact, no
non-trivial harmonic elements on $\widehat{SU(2)}$. Another
problem is that in general one expects a boundary of a
non-commutative space to be non-commutative. So if one wants to
generalize Biane's result, one needs not only to describe the pure
states of the algebra corresponding to extremal harmonic elements,
but also to justify its multiplicative structure.

One of the main properties of the Martin boundary is that any
harmonic element can be represented as an integral of the Martin
kernel by some measure. The unit can be represented by a canonical
measure, and the Martin boundary, considered as a measure space,
is called the Poisson boundary, see e.g. \cite{KV,K2}. It turns
out that the algebra of bounded measurable functions on the
Poisson boundary is canonically isomorphic to the space of bounded
harmonic elements, which a priori is just an operator space, but
actually has a unique structure of a von Neumann algebra. This
allows to define the Poisson boundary in the non-commutative
setting, namely by postulating that the algebra of bounded
measurable functions on it is the algebra of bounded harmonic
elements. This approach was suggested by Izumi \cite{I}, who was
motivated by subfactor theory and ITP-actions of compact quantum
groups. As we already mentioned, there are no non-constant
harmonic elements in $L(G)$ \cite{B1,B2,I}, so the Poisson
boundary of the dual of a compact group is trivial. The situation
changes drastically if we remove the assumption of cocommutativity
and instead of $L(G)$ consider the algebra of bounded functions on
a discrete quantum group. Izumi proved that the Poisson boundary
is non-trivial for any non-Kac algebra. One of the principal
results in \cite{I} is the computation of the Poisson boundary of
$\widehat{SU_q(2)}$. For this purpose Izumi carried out a detailed
study of the Markov operator $P_\phi$ for the $q$-trace $\phi$
associated with the fundamental corepresentation of $SU_q(2)$.
This allowed him to prove that the Poisson boundary is (the weak
operator closure of) the quantum homogeneous sphere of Podle\'{s}.
Then he extended this result to $q$-traces with finite support.

The appearance of $2$-spheres in the works of Biane and Izumi is
of course no coincidence. However, their interpretations of the
spheres as boundaries are quite different, and in the case covered
by both authors ($\phi$ is a state and $q=1$) the boundaries are
trivial. The missing link is the theory of the Martin boundary for
discrete quantum groups. The main objective in this paper is to
develop such a theory.

\smallskip

The paper is organized as follows.

In Section \ref{s1} we gather some facts about quantum groups,
with proofs included in cases we could not find a good reference.
Here we also discuss Radon-Nikodym cocycles for actions of
discrete quantum groups. They will play an important role in our
considerations for the same reason as in the classical theory,
where the Poisson kernel can be described as the Radon-Nikodym
cocycle of a measure representing the unit.

In Section \ref{s2} we study the Markov operator $P_\phi$ given by
convolution with a $q$-tracial state~$\phi$. We develop the part of the
potential theory needed to construct the boundary, featuring the
balayage theorem, which gives a canonical approximation of
harmonic elements by potentials $\sum_{n=0}^\infty P^n_\phi (x)$.
Although the adaptation of this theorem to our setting is fairly
straightforward, it is nevertheless striking that the lattice
property for selfadjoint elements is not needed. Note that to talk
about potentials we need to assume transience, i.e., that the
expected number of returns of the random walk to the origin is
finite. This condition is also necessary for the existence of
non-constant harmonic elements. It turns out that it is fulfilled
automatically in the generic quantum group case. Namely, the
Markov operators we consider always have a positive eigenvector
given by quotients of the classical and the quantum dimensions
with eigenvalue strictly less than one in the generic case. This
does not only imply transience, but also shows that the
probability of return to the origin at the $n$th step decreases
exponentially. In the last subsection of Section \ref{s2} we
briefly review the results of Izumi for the Poisson boundary with
emphasis on the quantum path spaces of our random walks and on the
non-commutative $0$-$2$ law rather than on fusion algebras and
ITP-actions.

In Section \ref{s3} we define the Martin kernel $K$. It is a
completely positive map from the algebra of finitely supported
functions to the algebra of bounded functions on the discrete
quantum group. Then we define the Martin compactification
$\tilde{A}_\phi$ as the C$^*$-algebra generated by the algebra
$\hat A$ of functions vanishing at infinity and the image of $K$.
The Martin boundary $A_\phi$ is then the quotient C$^*$-algebra
$\tilde{A}_\phi/\hat A$. We prove that any harmonic element gives
rise to a positive linear functional on~$A_\phi$. The algebra
$A_\phi$ has a canonical time evolution, and there always exists a
state which represents the unit and has the KMS property with
respect to this evolution. Any state $\nu$ representing the unit
gives rise to a normal u.c.p. map $K^*$ from the von Neumann
algebra $\pi_\nu (A_\phi)''$ to the Poisson boundary.
We give sufficient conditions for
$K^*$ to be an isomorphism which respects the canonical actions of
the quantum group and its dual.

In Section \ref{s4} we consider the case of $SU_q(2)$. Under a
certain summability assumption on~$\phi$, we prove that the Martin
boundary $A_\phi$ together with the two actions of the quantum
groups $SU_q(2)$ and $\widehat{SU_q(2)}$ on it can be identified
with the quantum homogeneous sphere of Podle\'{s}. The proof is
inspired by Biane's argument for ordinary $SU(2)$ \cite{B4}. The
first step is to realize the sphere as the quotient by the
compacts of the algebra of invariant formal pseudodifferential
operators of order zero associated to the $4D_+$-calculus of
Woronowicz. Then using the classical theory for $\7Z$ we show that
the Martin boundary of the center consists of one point. Combining
this with detailed knowledge of the representation theory of
$SU_q(2)$, including the Clebsch-Gordan coefficients, we compute
explicitly the Martin kernel on certain generating elements. The
main new problem in our context is to identify the actions of
$\widehat{SU_q(2)}$. This is solved by studying the Radon-Nikodym
cocycles for the actions. As a corollary we recover (and, in fact,
extend to states with infinite support) the results of Izumi.

\bigskip\bigskip

\section{Quantum Groups} \label{s1}

In Subsections \ref{s1.1}--\ref{s1.4} we collect a number of
definitions and results on compact and discrete quantum  groups
stated in a unified notation. For more details see
\cite{woro4,woro1,woro3,VD1,VD2,KD,Po,E}.

In the following $M(A)$ will denote the multiplier algebra of a
C$^*$-algebra $A$. We usually use the same symbol for a map and
its extension to the multiplier algebra. The C$^*$-algebra of
compact operators on a Hilbert space $H$ is denoted by $B_0(H)$.
We shall use $\odot$ to distinguish the algebraic tensor product
from the spatial tensor product $\otimes$, but we always use
$\otimes$ to denote tensor products of maps, then understood
according to context.

\subsection{Compact Quantum Groups} \label{s1.1}

A compact quantum group $(A ,\D)$ is a unital C$^*$-algebra $A$
and a unital $*$-homomorphism $\D \colon A\rightarrow A\otimes A$
such that $(\D\otimes\iota)\D=(\iota\otimes\D)\D$ and that both
$\D(A)(A\otimes1)$ and $\D(A)(1\otimes A)$ are dense in $A\otimes
A$. To any compact quantum group there exists a unique state
$\varphi$ of $A$ which is left- and right-invariant, i.e.,
$(\iota\otimes \varphi)\D =\varphi(\cdot)1$ and
$(\varphi\otimes\iota)\D=\varphi(\cdot)1$, respectively. This
functional is called the Haar functional and is always assumed to
be faithful, so we are considering a reduced compact quantum group
in the sense of \cite{ku-u}.

A unitary corepresentation $U$ of $(A,\D)$ on a Hilbert space
$H_U$ is a unitary element of $M(A\otimes B_0(H_U))$ such that
$(\D\otimes\iota)(U)=U_{13}U_{23}$. Here we have used the
leg-numbering convention. Let $\A$ consist of the elements
$(\iota\otimes\omega)(U)$, where $U$ is a finite dimensional
unitary corepresentation and $\omega$ is a linear functional of
$B(H_U)$. Then $\A$ is a dense unital $*$-subalgebra of $A$ such
that $\D (\A )\subset\A\odot\A$. Let $\A'$ denote the space of all
linear functionals on $\A$. There exist $\varepsilon\in\A'$ and a
linear map $S\colon\A\rightarrow\A$ such that $(\A,\D)$ is a Hopf
$*$-algebra with counit $\varepsilon$ and coinverse $S$, so $S(a^*
)=S^{-1}(a)^*$ for $a\in\A$.

Let $\omega\in\A'$ and $a\in\A$. Define $\omega*a ,a*\omega\in\A$
and $\bar\omega\in\A'$ by $\omega*a=(\iota\otimes\omega )\D(a)$,
$a*\omega =(\omega\otimes\iota )\D(a)$ and
$\bar\omega(a)=\overline{\omega(a^*)}$. There exists a family
$\{f_z\}_{z\in\7C}$ of unital, multiplicative functionals on $\A$
uniquely determined by the following properties:
\benu{7}{F1} $z\mapsto f_z (a)$ is an entire function of
exponential growth on the right-half plane for $a\in\A$.
\benu{7}{F2} $(f_z\otimes f_w)\D=f_{z+w}$ and $f_0=\varepsilon$.
\benu{7}{F3} $f_z S =f_{-z}$ and $\bar{f}_z =f_{-\bar{z}}$.
\benu{7}{F4} $S^2(a)=f_{-1}*a*f_1$ for $a\in\A$.
\benu{7}{F5} $\varphi(ab)=\varphi(b(f_1*a*f_1))$ for $a,b\in\A$.

The modular group $\{\sigma^\varphi_t\}_t$ of $\varphi$ is thus
determined by $\sigma^\varphi_t(a)=f_{it}*a*f_{it}$ for $a\in\A$
and $t\in\7R$. The one-parameter automorphism group $\{\tau_t\}_t$
of $A$ determined by $\tau_t(a)=f_{-it}*a*f_{it}$, for $a\in\A$
and $t\in\7R$, is called the scaling group of $(A,\D)$. The
involutive $*$-antiautomorphism $R$ of $A$ given by
$R(a)=f_\2*S(a)*f_{-\2}$ for $a\in\A$ is called the unitary
antipode. By definition $S=R\tau_{-\ii}=\tau_{-\ii}R$.

Let $U$ be a unitary corepresentation on $H_U$. Denote by $\bar
H_U$ the conjugate Hilbert space. Let~$J$ be the canonical
antilinear isometry $H_U\to\bar H_U$ and $j\colon B(H_U)\mapsto
B(\bar H_U)$, $j(x)=Jx^*J^{-1}$, the corresponding
$*$-antiisomorphism. Then the conjugate unitary corepresentation
$\bar U$ on $H_{\bar U}=\bar H_U$ is defined by $\bar U=(R\otimes
j)(U)$. The tensor product of two unitary corepresentations $U$
and $V$ is the unitary corepresentation $U\times V$ on $H_U\otimes
H_V$ defined by $U\times V=U_{12}V_{13}$.

Let $(H_\varphi,\xi_\varphi,\pi_r)$ be the GNS-triple associated
with the Haar functional $\varphi$. Since $\varphi$ is assumed to
be faithful we can identify $A$ with its image $\pi_r(A)$ in
$B(H_\varphi)$. The mapping $\xi\otimes a
\xi_\varphi\mapsto\D(a)(\xi\otimes\xi_\varphi)$, $\xi\in
H_\varphi$, $a\in A$, extends to a unitary operator on
$H_\varphi\otimes H_\varphi$. Its adjoint $W$ is called the
multiplicative unitary associated with $(A,\D)$. It is a unitary
corepresentation of $(A,\D)$ on $H_\varphi$ such that
$\D(a)=W^*(1\otimes a)W$ for $a\in A$. These two facts imply that
$W$ satisfies the pentagon equation
$W_{12}W_{13}W_{23}=W_{23}W_{12}$.

\smallskip

If $G$ is a compact group, then the C$^*$-algebra $C(G)$ of
continuous functions on $G$ with comultiplication $\D\colon
C(G)\to C(G)\otimes C(G)\cong C(G\times G)$ given by
$\D(f)(g,h)=f(gh)$ is a compact quantum group. If $\Gamma$ is a
discrete group, then the reduced group C$^*$-algebra
$C^*_r(\Gamma)$ of the group $\Gamma$ with comultiplication
$\D(\lambda_g)=\lambda_g\otimes\lambda_g$ is a compact quantum
group.

\subsection{Dual Discrete Quantum Groups} \label{s1.2}

Suppose $(A,\D)$ is a compact quantum group. Consider $\A'$ as a
unital $*$-algebra with product $\omega\eta
=(\omega\otimes\eta)\D$, unit $\varepsilon$ and $*$-operation
$\omega^* =\bar{\omega}S$. Define the $*$-subalgebra $\hat{\A}$ of
$\A'$ by $\hat{\A}=\{\hat{a}=a\varphi\ |\ a\in \A\}$, where
$(a\varphi)(b)=\varphi(ba)$ for $b\in\A$. For each
$\omega\in\hat\A$ the linear operator
$b\xi_\varphi\mapsto(b*\omega S^{-1})\xi_\varphi$ extends to a
bounded operator $\hat\pi_r (\omega )$ on $H_\varphi$ and
$\hat\pi_r$ is a faithful $*$-representation of $\hat\A$ on
$H_\varphi$. In the sequel we suppress the $*$-isomorphism
$\hat\pi_r$ of $\hat\A$ onto its image $\hat\pi_r(\hat\A)\subset
B(H_\varphi)$. Let $\hat A$ denote the norm closure of $\hat\A$ in
$B(H_\varphi)$. The formula $\hat\D (x)=W(x\otimes 1)W^*$ defines
a comultiplication $\hat\D\colon\hat A\to M(\hat A\otimes\hat A)$,
so $(\hat\D\otimes\iota )\hat\D=(\iota\otimes\hat\D)\hat\D$ and
both $\hat\D(\hat A)(\hat A\otimes1)$ and 
$\hat\D(\hat A)(1\otimes\hat A)$
are dense in $\hat A\otimes\hat A$. Henceforth $(\hat A,\hat\D)$
is called the discrete quantum group dual to $(A,\D)$. Note that
$W\in M(A\otimes\hat A)$. Any element $a\in\A$ extends to a
bounded linear functional on $\hat A$, so the latter can be
considered as a subspace of $\A'$. Moreover, $\hat A$ is a
$*$-subalgebra of $\A'$ and $\hat\D$ extends to a unital
$*$-homomorphism $\hat\D\colon\A'\to(\A\odot\A)'$ given by
$\hat\D(\omega)=\omega m$, where $m\colon\A\odot\A\rightarrow\A$
is the multiplication. The map $\hat\D$ satisfies the
coassociativity identity with counit and coinverse given by
$\hat\varepsilon (\omega)=\omega (1)$ and $\hat{S}(\omega)=\omega
S$, respectively. There exists a unique right-invariant Haar
weight $\hat\psi$ for $(\hat A,
\hat\D)$ which is a lower semicontinuous extension of the linear
functional on $\hat\A$ given by $\hat a\mapsto\varepsilon (a)$. It
satisfies the Plancherel formula $\hat\psi(\hat a^*\hat
b)=\varphi(a^* b)$ for $a,b\in\A$. Let $\rho=f_1$. It is a
positive self-adjoint operator affiliated with the von Neumann
algebra $\hat M$ generated by $\hat A\subset B(H_\phi)$. Since
$f_1$ is multiplicative, we have $\hat\D(\rho)=\rho\otimes\rho$
and $\hat S(\rho)=\rho^{-1}=f_{-1}$. Let $\hat\varphi
(x)=\hat\psi(\rho x\rho)$ for $x\in\hat A_+$, where $\hat A_+$ is
the cone of positive elements of the C$^*$-algebra $\hat A$. Then
$\hat\varphi$ is the (unique up to a scalar) left-invariant Haar
weight of $(\hat A,\hat\D)$. The modular group of $\hat\varphi$ is
given by $\sigma^{\hat\varphi}_t (x)=\rho^{it}x\rho^{-it}$ for
$x\in\hat A$, whereas that of $\hat\psi$ is given by
$\sigma^{\hat\psi}_t (x)=\rho^{-it}x\rho^{it}$. The scaling group
$\hat\tau$ coincides with $\sigma^{\hat\varphi}$, so the unitary
antipode $\hat R$ for $(\hat A,\hat\D)$ is given by $\hat
R(x)=\rho^{-\2}\hat S(x)\rho^\2$ for $x\in\hat A$.

Assume $U$ is a unitary corepresentation of $(A,\D)$. Then the
formula $\pi_U(\omega)=(\omega\otimes\iota)(U)$ for $\omega\in\hat
A$, defines a non-degenerate $*$-representation of $\hat
A\subset\A'$ on $H_U$. The map $U\mapsto\pi_U$ is an equivalence
of the category of unitary corepresentations of $(A,\D)$ and the
category of non-degenerate $*$-representations of $\hat A$. Note
that $\pi_W =\hat\pi_r$, and that the formula $\pi\mapsto
(\iota\otimes\pi)(W)$ defines an inverse functor. Moreover, these
functors preserve tensor products, where the tensor product of two
non-degenerate $*$-representations $\pi_1$ and $\pi_2$ is given by
$\pi_1\times\pi_2=(\pi_1\otimes\pi_2)\hat\D$, so $\pi_{U\times
V}=\pi_U\times\pi_V$.

\smallskip

Suppose $G$ is a compact group and $\Gamma$ is a discrete group.
The discrete quantum group dual to $(C(G),\D)$ is the group
C$^*$-algebra $C^*(G)$ of the group $G$ with comultiplication
$\hat\D(\int_Gf(g)\lambda_gdg)
 =\int_Gf(g)(\lambda_g\otimes\lambda_g)dg$
under the identification $\hat f=\int_Gf(g)\lambda_gdg$. The
discrete quantum group dual to $(C^*_r(\Gamma),\D)$ is the
C$^*$-algebra $c_0(\Gamma)$ of functions on $\Gamma$ vanishing at
infinity with comultiplication $\hat\Delta\colon c_0(\Gamma)\to
M(c_0(\Gamma)\otimes c_0(\Gamma)) =l^\infty(\Gamma\times\Gamma)$
given by $\hat\D(f)(g,h)=f(gh)$ under the identification
$\hat\lambda_g=\delta_{g^{-1}}$.

\subsection{Matrix Units and the Fourier Transform} \label{s1.3}

Denote by $I$ the set of equivalence classes of irreducible (and
thus finite dimensional) unitary corepresentations of a compact
quantum group $(A,\D)$, and let $U^s\in A\otimes B(H_s)$ denote a
fixed representative for the equivalence class $s\in I$. The
corresponding irreducible $*$-representation of $\hat A$ is
denoted by $\pi_s$. Then $\oplus_s \pi_s$ is a $*$-isomorphism of
$\hat A$ and the C$^*$-algebraic direct sum $\oplus_s B(H_s )$.
Under this isomorphism $\hat\A$ is the algebraic direct sum of
$B(H_s)$, $s\in I$, and $\A'$ is the algebraic direct product of
$B(H_s)$, $s\in I$. Denote by $I_s$ the unit of $B(H_s)$
considered as an element of $\hat\A$, so $\pi_s(x)=xI_s=I_sx$ for
$x\in\A'$. Then $\hat\psi(x)=\sum_sd_s\Tr\,\pi_s(x\rho^{-1})$,
where $d_s=\Tr\,\pi_s(\rho)=\Tr\,\pi_s(\rho^{-1})$ is the quantum
dimension of $U^s$ and $\Tr$ is the canonical, non-normalized
trace on $B(H_s)$. For each $s\in I$, we pick an orthonormal basis
$\{\xi^s_i\}_i$ for $H_s$. Let $m^s_{ij}$ denote the corresponding
matrix units in $B(H_s)$ with respect to $\{\xi^s_i\}_i$, so
$m^s_{ij}\xi^s_k=\delta_{jk}\xi^s_i$ and therefore
$m^s_{ij}m^s_{kl}=\delta_{jk}m^s_{il}$ and
${m^s_{ij}}^*=m^s_{ji}$. Let $u^s_{ij}\in\A$ be the matrix
coefficients of $U^s$, so $U^s =\sum_{i,j}u^s_{ij}\otimes
m^s_{ij}$. The fact that $U^s$ is a corepresentation means then
that $\D(u^s_{ij})=\sum_ku^s_{ik}\otimes u^s_{kj}$. Note that $U^s
=(\iota\otimes\pi_s)W$ implies
\begin{equation}
\label{e1.1}
W=\sum_s\sum_{i,j}u^s_{ij}\otimes m^s_{ij}.
\end{equation}
The following
orthogonality relations hold:

\begin{equation} \label{e1.2}
\varphi(({u^t_{kl}})^*u^s_{ij})
=\delta_{st}\delta_{jl}{f_{-1}(u^s_{ik})\over d_s}.
\end{equation}

To simplify some formulas we pick the basis $\{\xi^s_i\}_i
$ such that the matrix
$\pi_s (\rho )$ is diagonal. Then $f_z (u^s_{ij})=0$ if $i\ne j$,
and $\pi_s (\rho^z)=\sum_i f_z (u^s_{ii})m^s_{ii}$. We also have
\begin{equation} \label{e1.3}
\tau_t(u^s_{kl})=f_{it}(u^s_{kk})f_{-it}(u^s_{ll})u^s_{kl}\ \ \
\hbox{and}\ \ \
\sigma^\varphi_t(u^s_{kl})=f_{it}(u^s_{kk})f_{it}(u^s_{ll})u^s_{kl}.
\end{equation}

Let $\bar s$ be the equivalence class of $\overline{U^s}$. Since
$\overline{U^s}$ is equivalent to $U^{\bar s}$, there exists an
antilinear isometry $J_s\colon H_s\to H_{\bar s}$ uniquely defined
up to a scalar of modulus one. So $j_s(x)=J_sx^*J_s^{-1}$ is a
well-defined $*$-antiisomorphism from $B(H_s)$ to $B(H_{\bar s})$.
Let $n^{\bar s}_{ij}=j_s(m^s_{ji})$ be the new system of matrix
units in $B(H_{\bar s})$ and $v^{\bar s}_{ij}$ be the
corresponding matrix coefficients of $U^{\bar s}$, so $U^{\bar
s}=\sum_{i,j}v^{\bar s}_{ij}\otimes n^{\bar s}_{ij}$.

The Fourier transform $\F\colon\A\to\hat\A$ is the bijection
defined by $\F(a)=\hat a$.

\begin{lem}
\label{1.1}
The following formulas hold:
\enu{i} $m^s_{ij}(u^t_{kl})=\delta_{st}\delta_{ik}\delta_{jl}$;
\vspace{2mm}
\enu{ii} $v^{\bar
s}_{ij}=f_\2(u^s_{ii})f_{-\2}(u^s_{jj}){u^s_{ij}}^*$;
\vspace{2mm}
\enu{iii}
$\F(u^s_{ij})=d_s^{-1}f_{-\2}(u^s_{ii})f_{-\2}(u^s_{jj})n^{\bar
s}_{ij}$;
\vspace{2mm}
\enu{iv}
$\hat R(x)=j_s(x)$ for $x\in B(H_s)$, so $\hat
S(m^s_{ij})=f_\2(u^s_{ii})f_{-\2}(u^s_{jj})n^{\bar s}_{ji}$.
\end{lem}

\bp
Part (i) is obvious, because by definition we have
$\delta_{st}m^s_{ij}=\pi_t(m^s_{ij})
=\sum_{k,l}m^s_{ij}(u^t_{kl})m^t_{kl}$.

By definition of $\overline{U^s}$ we have $v^{\bar
s}_{ij}=R(u^s_{ji})$. Since $S(u^s_{ji})={u^s_{ij}}^*$,
$S=R\tau_{-\ii}$ and
$\tau_\ii(u^s_{ji})=f_{-\2}(u^s_{jj})f_\2(u^s_{ii})u^s_{ji}$ by
(\ref{e1.3}), we get (ii).

By definition of the Fourier transform the orthogonality relations
can be rewritten as
$$
\F(u^s_{ij})({u^t_{kl}}^*)
=\delta_{st}\delta_{jl}\delta_{ik}f_{-1}(u^s_{ii})d_s^{-1},
$$
or, in view of (ii), as $\F(u^s_{ij})(v^{\bar
t}_{kl})=\delta_{st}\delta_{jl}\delta_{ik}f_{-\2}(u^s_{ii})
f_{-\2}(u^s_{jj})d_s^{-1}$. Since $n^{\bar s}_{ij}(v^{\bar
t}_{kl})=\delta_{st}\delta_{ik}\delta_{jl}$ by (i), we get (iii).

By definition of the conjugate unitary corepresentation we have
$$
\pi_{\bar U}(\omega)\bar\xi=\overline{\pi_U(\omega R)^*\xi}
=\overline{\pi_U(\hat R(\omega^*))\xi}\ \ \ \hbox{for}\
\xi\in H_U\ \hbox{and}\ \omega\in\hat A.
$$
This implies the first part of (iv). The second part follows from
$\hat S=\hat R(\rho^\2\cdot\rho^{-\2})$.
\ep

We denote by $0$ the equivalence class of the trivial
one-dimensional corepresentation, so $\pi_0=\hat\varepsilon$. Note
that $\varphi=\F (1)=I_0$ because $I_0^2=I_0$ and $\hat a
I_0=\hat\varepsilon(\hat a)I_0$, for $a\in\A$, determine $I_0$
uniquely.

\begin{lem} \label{1.2}
We have
$$
\hat\D(I_0)=\sum_s\sum_{i,j}\F(u^s_{ij})\otimes m^s_{ij}
=\sum_s\sum_{i,j}d_s^{-1}f_{-\2}(u^s_{ii})
f_{-\2}(u^s_{jj})n^{\bar s}_{ij}\otimes m^s_{ij}.
$$
\end{lem}

\bp
This is easily verified by applying the functionals on both sides
of the identity to the linear basis $\{v^s_{ij}\otimes
u^t_{kl}\}_{s,t,i,j,k,l}$ of $\A\odot\A$ and using the
orthogonality relations (\ref{e1.2}) together with
Lemma~\ref{1.1}. Alternatively, as $(\iota\otimes b)(W)=b$ for any
$b\in\A$, we have $(a\otimes
b)(\F\otimes\iota)(W)=a\F(b)=\varphi(ab)$ and $(a\otimes
b)\hat\D(I_0)=I_0(ab)=\varphi(ab)$ for all $a,b\in\A$. Hence
$\hat\D(I_0)=(\F\otimes\iota)(W)$.
\epp

\subsection{Coactions and Invariant States} \label{s1.4}

Suppose $(A,\D)$ is a compact or discrete quantum group. A left
(resp. right) coaction $\alpha$ of $(A,\D)$ on a C$^*$-algebra $B$
is a $*$-homomorphism $\alpha\colon B\to M(A\otimes B)$ such that
$\alpha (B)(A\otimes 1)$ (resp. $\alpha (B)(1\otimes A)$) is dense
in $A\otimes B$ (resp. $B\otimes A$) and that
$(\iota\otimes\alpha)\alpha =(\D\otimes\iota)\alpha$ (resp.
$(\alpha\otimes\iota)\alpha =(\iota\otimes\D)\alpha$).

The fixed point algebra $B^\alpha$ for a left coaction $\alpha$ is
the C$^*$-subalgebra of $B$ consisting of elements $x\in B$ such
that $\alpha(x)=1\otimes x$.

\begin{prop} \label{1.3}
Let $\alpha\colon B\to M(B\otimes A)$ be a right coaction of
$(A,\D)$ on $B$. Define $\B=B$ in the discrete case and
$\B=\span\{(\iota\otimes\varphi)(\alpha (b)(1\otimes a))\ |\ a\in\A,\
b\in B\}$ in the compact case. Then $\B$ is a dense $*$-subalgebra
of $B$ and
$\B\odot\A=\alpha(\B)(1\otimes\A)=(1\otimes\A)\alpha(\B)$.
Moreover, the operator $r\colon\B\odot\A\to\B\odot\A$, $r(b\otimes
a)=\alpha(b)(1\otimes a)$, is invertible with inverse given by
$s(b\otimes a)=(\iota\otimes S)((1\otimes S^{-1}(a))\alpha(b))$.
\end{prop}

\bp
In the compact case the orthogonality relations (\ref{e1.2}) imply
(see e.g. \cite[Theorem~1.5]{Po}) that $\B$ is dense in~$B$, 
$\alpha (\B)\subset\B\odot\A$ and
$$
\B=\span\{(\iota\otimes\varphi)(\alpha (b)(1\otimes a))\
|\ a\in\A,\ b\in\B\}.
$$
In the discrete case $(B\otimes A)(1\otimes I_s)=B\odot
B(H_s )$, so $\alpha (\B)(1\otimes\A)\subset\B\odot\A$ and
$(1\otimes\A)\alpha (\B)\subset\B\odot\A$.

Let us prove that $(\iota\otimes\varepsilon)\alpha (b)=b$ for
$b\in\B$. Note that in the discrete case $\varepsilon$ is defined
on $A$, whereas in the compact case the formula makes sense due
to the property $\alpha(\B)\subset\B\odot\A$.
Apply $\iota\otimes\varepsilon\otimes\iota$ to
the identity
$$
(\iota\otimes\D)(\alpha (b)(1\otimes a))
=(\alpha\otimes\iota)\alpha(b)(1\otimes\D(a)),
$$
which yields
$$
\alpha(b)(1\otimes a)
=(\iota\otimes\varepsilon\otimes\iota)
((\alpha\otimes\iota)\alpha(b))(1\otimes a)
=((\iota\otimes\varepsilon)\alpha\otimes\iota)
(\alpha(b)(1\otimes a)).
$$
In the discrete case the density of $\alpha(\B)(1\otimes\A)$ in
$B\otimes A$ gives the result, whereas in the compact case we
apply $\iota\otimes\varphi$ and use the description of $\B$ given
above.

Next consider $a\in\A$ and $b\in\B$. Observe that $r=(\iota\otimes
m)(\alpha\otimes\iota)$. Thus since
$(\iota\otimes\varepsilon)\alpha(b)=b$ and $m(\iota\otimes
S)\D=\varepsilon(\cdot)1$ we get
\begin{eqnarray*}
rs(b\otimes a)
&=&(\iota\otimes m)(\alpha\otimes\iota)
(\iota\otimes S)((1\otimes S^{-1}(a))\alpha(b))\\
&=&(\iota\otimes m)(\iota\otimes\iota\otimes S)
    ((\iota\otimes\iota\otimes S^{-1}(a))
    (\alpha\otimes\iota)\alpha(b))\\
&=&(\iota\otimes m)(\iota\otimes\iota\otimes S)
((\iota\otimes\iota\otimes S^{-1}(a))(\iota\otimes\D)\alpha(b))\\
&=&(\iota\otimes\varepsilon(\cdot)1)\alpha(b)(1\otimes a)
=b\otimes a,
\end{eqnarray*}
so $rs=\iota$. The equality $sr=\iota$ is proved analogously.
\ep

The proof of the previous proposition shows the following.

\begin{cor}
\label{1.4}
\mbox{\ }
\enu{i}
Let $(A,\D)$ be a compact quantum group and $\alpha\colon B\to
B\otimes A$ a $*$-homomorphism such that
$(\alpha\otimes\iota)\alpha =(\iota\otimes\D)\alpha$. Then
$\alpha$ is a coaction if and only if there exists a dense
$*$-subalgebra $\B$ of $B$ such that $\alpha(\B)\subset\B\odot\A$
and $(\iota\otimes\varepsilon)\alpha=\iota$ on $\B$.
\enu{ii}
Let $(A,\D)$ be a discrete quantum group and $\alpha\colon B\to
M(B\otimes A)$ a non-degenerate $*$-homomor\-phism such that
$(\alpha\otimes\iota)\alpha =(\iota\otimes\D)\alpha$. Then
$\alpha$ is a coaction if and only if
$(\iota\otimes\varepsilon)\alpha=\iota$.
\epp
\end{cor}

An invariant state for a right coaction $\alpha$ of a compact or
discrete quantum group $(A,\D)$ on $B$ is a state $\eta$ on $B$
such that $(\eta\otimes\iota)\alpha=\eta(\cdot)1$.

\begin{prop}
\label{1.5}
Let $\eta$ be an invariant state for a right coaction $\alpha$ of
a compact or discrete quantum group. Then
\enu{i}
$(\eta\otimes\omega)(\alpha (b_1)(b_2\otimes 1))
=(\eta\otimes\omega S)((b_1\otimes 1)\alpha(b_2))$ for any bounded
functional $\omega$ of $A$ such that the functional $\omega S$ on
$\A$ extends to a bounded functional on $A$;
\enu{ii}
$\alpha\sigma^\eta_t =(\sigma^\eta_t\otimes\tau_{-t})\alpha$
whenever $\eta$ is a faithful KMS-state (at inverse temperature
$\beta=-1$).
\end{prop}

\bp
Let $(H_\eta ,\xi_\eta, \pi_\eta)$ be a GNS-triple for $\eta$, and
as usual we suppress $\pi_\eta$. It is easy to check
that the formula
$$U(b\xi_\eta\otimes\xi)=\alpha(b)(\xi_\eta\otimes\xi),$$
for $b\in B$, $\xi\in H_\varphi$, defines a  unitary
corepresentation $U$ of $(A,\D)$ on $H_\eta$ such that $\alpha
(b)=U(b\otimes 1)U^*$. Then $(\iota\otimes\omega
S)(U)=(\iota\otimes\omega)(U^*)$ for example 
by \cite[Proposition~5.2]{ku-u}. Since
any bounded linear functional on $A$ can be weakly approximated by
vector functionals, by definition of $U$ we conclude that
$$
(\iota\otimes\omega)(U)b\xi_\eta
=(\iota\otimes\omega)\alpha (b)\xi_\eta .
$$
Hence
\begin{eqnarray*}
\eta((\iota\otimes\omega)\alpha(b_1)b_2)
&=&(b_2\xi_\eta ,(\iota\otimes\bar\omega)\alpha(b_1^*)\xi_\eta)
=(b_2\xi_\eta, (\iota\otimes\bar\omega)(U)b_1^*\xi_\eta)\\
&=&((\iota\otimes\omega)(U^*)b_2\xi_\eta ,b_1^*\xi_\eta)
=((\iota\otimes\omega S)(U)b_2\xi_\eta ,b_1^*\xi_\eta))\\
&=&((\iota\otimes\omega S)\alpha (b_2)\xi_\eta ,b_1^*\xi_\eta))
=(\eta\otimes\omega S)((b_1\otimes 1)\alpha(b_2)).
\end{eqnarray*}
This proves (i). Then (ii) follows for instance from the proof of
\cite[Theorem 2.9]{E}.
\ep

If $\eta$ is an invariant state for a left coaction $\alpha$, we
have the following analogous formulas:
\enu{i}
$(\omega\otimes\eta)(\alpha (b_1)(1\otimes b_2))
=(\omega S^{-1}\otimes\eta)((1\otimes b_1)\alpha(b_2))$,
\enu{ii}
$\alpha\sigma^\eta_t =(\tau_t\otimes\sigma^\eta_t)\alpha$ whenever
$\eta$ is a faithful KMS-state.

\medskip

Let $\eta$ be a KMS-state on a C$^*$-algebra $B$. The formula
$$
(x,y)_\eta=\eta(x\sigma^\eta_{-\ii}(y^*))
 =(xJ_\eta y\xi_\eta,\xi_\eta),
$$
where $J_\eta$ is the modular involution, defines an inner product
$(\cdot,\cdot)_\eta$ on the von Neumann algebra $N=\pi_\eta
(B)''$. If we are given two C$^*$-algebras $B_1$ and $B_2$ with
KMS-states $\eta_1$ and $\eta_2$, and a completely positive
contraction $T \colon B_1\to B_2$ such that $\eta_2 T=\eta_1$,
then there exists a unique normal unital completely positive map
$T^*
\colon N_2\to N_1$ such that $\eta_1 T^*=\eta_2$ and
$(x,T^*y)_{\eta_1}=(Tx,y)_{\eta_2}$ for $x\in B_1$ and $y\in B_2$
\cite{AC}.

\begin{prop} \label{1.6}
Suppose $\eta$ is a faithful invariant KMS-state for a right
coaction $\alpha$ of a discrete or compact quantum group $(A,\D)$
on a C$^*$-algebra $B$. Define $T_\omega  \colon B\to B$ by
$T_\omega =(\iota\otimes\omega)\alpha$, where~$\omega$ is a state
on $A$. Then $T_\omega$ is a completely positive contraction,
$\eta T_\omega=\eta$ and $T_\omega^* =T_{\omega R}$.
\end{prop}

\bp
By density of analytic elements it is sufficient to consider $x$
and $y$ to be $\sigma^\eta$-analytic and $\omega$ to be
$\tau$-analytic. Then Proposition \ref{1.5}(i) says that
$$
\eta(T_\omega (x)\sigma^\eta_t(y^*))
=\eta(xT_{\omega S}(\sigma^\eta_t(y^*))),
$$
while part (ii) of that proposition implies
$T_{\omega S}\sigma^\eta_t =\sigma^\eta_t T_{\omega S\tau_{-t}}$.
Thus
$$
\eta(T_\omega (x)\sigma^\eta_t(y^*))
=\eta(x\sigma^\eta_t T_{\omega S\tau_{-t}}(y^*))
=\eta(\sigma^\eta_{-t}(x)T_{\omega S\tau_{-t}}(y^*)).
$$
Taking analytic continuation to $t=-\ii$ and using $\omega
S\tau_{-t}=\omega\tau_{-t-\ii}R$ completes the proof.
\ep

{\it From this point onwards $(A,\D)$ will always denote a compact
quantum group}.

\medskip

We may clearly regard $\D$ as a left coaction of $(A, \D)$ on $A$
and $\hat\D$ as a right coaction of $(\hat A,
\hat\D)$ on $\hat A$. We shall consider two more coactions.

\begin{prop} \label{1.7}
The formulas
$$
\Phi (x)=W^* (1\otimes x)W\ \ \ \hbox{and}\ \ \
\hat\Phi (a)=W(a\otimes 1)W^*,$$
for $x\in\hat A$ and $a\in A$, define a left coaction $\Phi$ of
$(A,\D)$ on $\hat A$ and a right coaction $\hat\Phi$ of $(\hat
A,\hat\D)$ on~$A$.
\end{prop}

\bp
Coassociativity follows from the pentagon equation for $W$. Note
that $\Phi (x)={U^s}^*(1\otimes x)U^s$ for $x\in B(H_s)$, so the
result follows from Corollary \ref{1.4} and the properties
$(\varepsilon\otimes\iota)U^s=I_s$ and
$(\iota\otimes\hat\varepsilon)W=1$.
\ep

The coactions $\Phi$ and $\hat\Phi$ are analogues of the adjoint
action of a group on its group C$^*$-algebra. Suppose again that
$G$ is a compact group and $\Gamma$ is a discrete group. If
$(A,\D)=(C(G),\D)$, then $\hat\Phi$ is trivial,
$\hat\Phi(a)=a\otimes1$, while $\Phi\colon C^*(G)\to M(C(G)\otimes
C^*(G))\subset L^\infty(G,W^*(G))$ is given by
$\Phi(\int_Gf(g)\lambda_gdg)(h)=\int_Gf(g)\lambda_{h^{-1}gh}dg$.
On the other hand, if $(A,\D)=(C^*_r(\Gamma),\D)$, then $\Phi$ is
trivial, $\Phi(x)=1\otimes x$, while $\hat\Phi\colon
C^*_r(\Gamma)\to M(C^*_r(\Gamma)\otimes
c_0(\Gamma))=l^\infty(\Gamma,C^*_r(\Gamma))$ is given by
$\hat\Phi(\lambda_g)(h)=\lambda_{hgh^{-1}}$.

\smallskip

Consider the left coaction $\alpha_s =\Phi|_{B(H_s)}$
of $(A,\D)$ on $B(H_s)$. The restriction of the $q$-trace
$\phi_s=d_s^{-1}\Tr\,\pi_s(\cdot\rho^{-1})$ to $B(H_s)$ is an
invariant state for this coaction. To see this, first note that
$$
(S^2\otimes\iota)U^s
=(1\otimes\pi_s(\rho))U^s(1\otimes\pi_s(\rho^{-1}))
$$
by property (F4) for the family $\{f_z\}$, so since
$(S\otimes\iota)U^s ={U^s}^*$, we get
$$
(S^{-1}\otimes\iota)(U^s)=(S^{-2}\otimes\iota)({U^s}^*)
=(S^2\otimes\iota)(U^s)^*
=(1\otimes\pi_s(\rho^{-1})){U^s}^*(1\otimes\pi_s (\rho)).
$$
Thus
\begin{eqnarray*}
d_s(\iota\otimes\phi_s)\alpha_s(x)
&=&(\iota\otimes\Tr)({U^s}^*(1\otimes x)U^s
    (1\otimes\pi_s (\rho^{-1})))\\
&=&(\iota\otimes\Tr)((1\otimes x)(S\otimes\iota)
    ((S^{-1}\otimes\iota)(U^s)(1\otimes\pi_s(\rho^{-1}))
      (S^{-1}\otimes\iota)({U^s}^*)))\\
&=&(\iota\otimes\Tr)((1\otimes x)(S\otimes\iota)
     ((1\otimes\pi_s(\rho^{-1})){U^s}^*U^s ))=d_s\phi_s(x)1.
\end{eqnarray*}

Note that since $U^s$ is irreducible, we have
$B(H_s)^{\alpha_s}=\7C I_s$. The formula $E_s
=(\varphi\otimes\iota)\alpha_s$ defines a conditional expectation
of $B(H_s)$ onto $B(H_s)^{\alpha_s}$, and any invariant state
$\eta$ must satisfy the property $\eta E_s=\eta$. As a consequence
$E_s =\phi_s (\cdot)I_s$ and $\phi_s|_{B(H_s)}$ is the unique
invariant state for $\alpha_s$.

Let $\C\subset\hat{A}^*$ denote the norm closure of the linear
span of $\{\phi_s\}_{s\in I}$, so $\C$ consists of the functionals
$\sum_s\lambda_s\phi_s$, $\{\lambda_s\}_s\in l^1(I)$. Any state in
$\C$ is an invariant state for $\Phi$. Conversely, since
$\phi_s|_{B(H_s)}$ is the only invariant state for $\alpha_s$, any
invariant state for $\Phi$ belongs to $\C$.

For any unitary corepresentation $U$ of $(A,\D)$, the formula
$\alpha_U (x)=U^*(1\otimes x)U$, for $x\in B_0 (H_U)$, defines a
left coaction $\alpha_U$ of $(A,\D)$ on $B_0(H_U)$, so $\alpha_s
=\alpha_{U^s}$. Again by uniqueness of invariant states for
$\alpha_s$, for any invariant state $\eta$ on $B_0(H_U)$ we have
$\tilde{\eta}\pi_U\in\C$, where $\tilde\eta$ is the unique normal
extension of $\eta$ to $B(H_U)$.

\begin{lem} \label{1.8}
\mbox{\ }
\enu{i} The state $\phi_s$ considered as a linear functional on
$\hat\A$ lies in $\A$ and
$\phi_s=\sum_id_s^{-1}f_{-1}(u^s_{ii})u^s_{ii}$.
\enu{ii} The linear space $\C$ is a
subalgebra of ${\hat A}^*$. Namely, if $U^s\times
U^t\cong\sum_wN^w_{s,t}U^w$ is the decomposition of the
corepresentation $U_s\times U_t$ into irreducible ones, then
$\displaystyle\phi_s\phi_t=\sum_w{d_w\over
d_sd_t}N^w_{s,t}\phi_w$.
\end{lem}

\bp Statement (i) follows from Lemma~\ref{1.1}(i).
Since $\hat\D(\rho)=\rho\otimes\rho$, property (ii) follows from
$\sum_wN^w_{s,t}\Tr\,\pi_w=\Tr\,\pi_{U^s\times
U^t}=(\Tr\otimes\Tr)(\pi_s\otimes \pi_t)\hat\D$.
\ep

Note that in general $\C$ is not a $*$-subalgebra because
$\phi_s^*=d_{\bar s}^{-1}\Tr\,\pi_{\bar s}(\cdot\rho)$.

\subsection{Radon-Nikodym Cocycle}

Let $\alpha \colon B\to M(B\otimes\hat A)$ be a right coaction of
a discrete quantum group $(\hat A ,\hat\D)$ on a unital
C$^*$-algebra~$B$. Denote by $M(B\odot\hat\A)$ the algebraic
multiplier algebra of $B\odot\hat\A$ \cite{VD2}, so
$$
M(B\odot\hat\A)=\prod_{s\in I}B\otimes B(H_s),
$$
where we use the algebraic direct product. Thus $M(B\odot\hat\A)$
is a $*$-algebra, but not a normed algebra.

\begin{defin} \label{1.9}
A state $\eta$ on $B$ is called quasi-invariant (with respect to a
right coaction $\alpha \colon B\to M(B\otimes\hat A)$) with
Radon-Nikodym cocycle $y\in M(B\odot\hat\A)$ if
\enu{i}
$(\eta\otimes\iota)\alpha (b)=(\eta\otimes\iota)
((b\otimes 1)(\iota\otimes\hat S)(y))$ for $b\in B$;
\enu{ii}
$(\iota\otimes\hat\D)(y)=(\alpha\otimes\iota)(y)(y\otimes 1)$.
\end{defin}
Note that these formulas make sense if we extend all appropriate
homomorphisms to algebraic multiplier algebras. A more concrete
way to proceed is as follows. First observe that
$M(\hat\A)=\A'=\prod_{s\in I}B(H_s)$, so both $\hat S\colon
\A'\to\A'$ and the homomorphism $\hat\D\colon \A'\to
(\A\odot\A)'=\prod_{s,t\in I}B(H_s)\otimes B(H_t)$ are
well-defined, see Subsection~\ref{s1.2}. Since $\hat
S(B(H_s))=B(H_{\bar s})$, clearly $\iota\otimes\hat S$ is
well-defined. Thus identity (i) makes sense. Concerning (ii),
notice that for fixed $s,t\in I$, there exists only finitely many
$w\in I$ such that $\hat\D (x)$, $x\in B(H_w)$, has a non-zero
component in $B(H_s)\otimes B(H_t)$, namely, those $w\in I$ for
which $U^w$ is a subcorepresentation of $U^s\times U^t$. Thus
$(\iota\otimes\hat\D)(y)$ is a well-defined element of
$M(B\odot\hat\A\odot\hat\A)=\prod_{s,t\in I}B\otimes B(H_s)\otimes
B(H_t)$.

Observe that whenever $\eta$ is faithful, condition (i) of this
definition determines $y$ uniquely.

\begin{prop} \label{1.10}
Suppose $\eta$ is a state on $B$ and $y\in
M(B\odot\hat\A)$ satisfies identity (i) in Definition~\ref{1.9}.
Then:
\enu{i}
if $\eta$ is faithful, the element $y$ satisfies identity
(ii) in Definition \ref{1.9}, so $\eta$ is quasi-invariant with
Radon-Nikodym cocycle $y$;
\enu{ii}
$(\eta\otimes\iota)(\alpha(b_1)(b_2\otimes 1))
=(\eta\otimes\hat S)((b_1\otimes 1)\alpha (b_2)y)$
for $b_1 ,b_2\in B$.
\end{prop}

\bp
We begin by proving (ii). Let $z=b\otimes a$ for $b\in B$ and
$a\in\hat\A$, and consider the linear maps $r,s$ on $B\odot\hat\A$
introduced in Proposition \ref{1.3}. Then
\begin{eqnarray*}
(\eta\otimes\iota)(\alpha(b_1)r(z))
&=&(\eta\otimes\iota)(\alpha(b_1 b)(1\otimes a))\\
&=&(\eta\otimes\iota)((b_1 b\otimes 1)(\iota\otimes\hat S)(y)
     (1\otimes a))\\
&=&(\eta\otimes\hat S)((b_1\otimes 1)(\iota\otimes\hat S^{-1})(z)y).
\end{eqnarray*}
If we apply this identity to $z=s(b_2\otimes I_t)$, and use
$rs=\iota$, we get
\begin{eqnarray*}
I_t (\eta\otimes\iota)(\alpha(b_1)(b_2\otimes 1))&=&
(\eta\otimes\iota)(\alpha(b_1)(b_2\otimes I_t))\\
&=&(\eta\otimes\hat S)((b_1\otimes 1)(\iota\otimes\hat S^{-1})
     s(b_2\otimes I_t)y)\\
&=&(\eta\otimes\hat S)((b_1\otimes 1)(1\otimes I_{\bar t})
     \alpha (b_2 )y)\\
&=&I_t(\eta\otimes\hat S)((b_1\otimes 1)\alpha (b_2 )y).
\end{eqnarray*}
Since this holds for any $t\in I$, assertion (ii) is proved.

\smallskip

To prove (i), apply $\hat\D$ to the identity
$$
(\eta\otimes\iota)\alpha (b)
=(\eta\otimes\iota)((b\otimes 1)(\iota\otimes\hat S)(y)).
$$
The right hand side yields $(\eta\otimes\chi)(\iota\otimes\hat
S\otimes\hat S) ((b\otimes 1\otimes 1)(\iota\otimes\hat\D)(y))$,
where $\chi$ denotes the flip on $(\A\odot\A)'$. Whereas the
left hand side gives
\begin{eqnarray*}
(\eta\otimes\iota\otimes\iota)(\iota\otimes\hat\D)\alpha(b)
&=&(\eta\otimes\iota\otimes\iota)(\alpha\otimes\iota)\alpha(b)\\
&=&(\eta\otimes\iota\otimes\iota)(\alpha(b)_{13}
      (\iota\otimes\hat S\otimes\iota)(y\otimes 1))\\
&=&(\eta\otimes\chi)((\alpha(b)\otimes 1)
      (\iota\otimes\iota\otimes\hat S)(y_{13}))\\
&\stackrel{(ii)}{=}&(\eta\otimes\chi)
     (\iota\otimes\hat S\otimes\iota)
       ((b\otimes 1\otimes 1)(\alpha\otimes\iota)
             (\iota\otimes\hat S)(y)(y\otimes 1))\\
&=&(\eta\otimes\chi)(\iota\otimes\hat S\otimes\hat S)
     ((b\otimes 1\otimes 1)(\alpha\otimes\iota)(y)(y\otimes 1)).
\end{eqnarray*}
Thus
$$
(\eta\otimes\chi)(\iota\otimes\hat S\otimes\hat S)
((b\otimes 1\otimes 1)(\iota\otimes\hat\D)(y))
=(\eta\otimes\chi)(\iota\otimes\hat S\otimes\hat S)
((b\otimes 1\otimes 1)(\alpha\otimes\iota)(y)(y\otimes 1)).
$$
Since $b$ is arbitrary and $\eta$ is faithful, assertion (i) now
follows.
\ep

We call property (ii) strong quasi-invariance, so (i) in
Proposition~\ref{1.5} follows from Proposition~\ref{1.10} (at
least in the discrete case).

In general $y$ is not self-adjoint. But it will be self-adjoint
and have other nice properties under additional assumptions which
will always be satisfied in our examples.

\begin{prop} \label{1.11}
Let $\eta$ be a quasi-invariant state with Radon-Nikodym cocycle
$y$. Suppose in addition that $\eta$ is a faithful KMS-state and
that $\alpha\sigma^\eta_t
=(\sigma^\eta_t\otimes\hat\tau_{-t})\alpha$. Then
\enu{i}
$(\sigma^\eta_t\otimes\hat\tau_{-t})(y)=y$;
\enu{ii}
$y$ is positive and invertible with $y^{-1}=(\iota\otimes
m)(\alpha\otimes\iota)(\iota\otimes\hat S)(y)$.
\end{prop}

\bp
Apply $\hat\tau_{-t}$ to the identity
$(\eta\otimes\iota)\alpha (b)=(\eta\otimes\iota)
((b\otimes 1)(\iota\otimes\hat S)(y))$,
and use $\sigma^\eta_t$-invariance of~$\eta$. Thus
$$
(\eta\otimes\iota)\alpha\sigma^\eta_t(b)
=(\eta\otimes\iota)((\sigma^\eta_t(b)\otimes 1)
(\sigma^\eta_t\otimes\hat\tau_{-t})(\iota\otimes\hat S)(y))
=(\eta\otimes\iota)((\sigma^\eta_t(b)\otimes 1)
(\iota\otimes\hat S)(\sigma^\eta_t\otimes\hat\tau_{-t})(y)),
$$
for $b\in B$, and (i) follows by uniqueness.

\smallskip

To see that $y\geq 0$, it is clearly sufficient to show that
$\sum_{k,j}(\eta\otimes\iota)((b_k^*\otimes a_k^*)y(b_j\otimes
a_j))\geq 0$ for $a_k\in\hat A$ and $b_k\in B$. Moreover, we can
assume that $b_k$ is $\sigma^\eta$-analytic. Write $x_k
=\sigma^\eta_\ii (b_k)$. Using (i) and the KMS-condition, we thus
get

\vspace{2mm}
$\displaystyle
\sum_{k,j}(\eta\otimes\iota)((b_k^*\otimes a_k^*)y(b_j\otimes a_j))$
\vspace{-4mm}
\begin{eqnarray*}
&=&\sum_{k,j}a_k^*(\eta\otimes\iota)
     ((\sigma^\eta_i(b_j)b_k^*\otimes 1)y )a_j
=\sum_{k,j}a_k^*(\eta\otimes\iota)
      ((x_j x_k^*\otimes 1)(\sigma^\eta_{-\ii}\otimes\iota)(y))a_j\\
&=&\sum_{k,j}a_k^*(\eta\otimes\iota)
      ((x_j x_k^*\otimes 1)(\iota\otimes\hat\tau_{-\ii})(y))a_j
=\sum_{k,j}a_k^*(\eta\otimes\hat R)
      ((x_j x_k^*\otimes 1)(\iota\otimes\hat S)(y))a_j\\
&=&\sum_{k,j}a_k^*(\eta\otimes\hat R)\alpha(x_j x_k^*)a_j
=\sum_{k,j}(\eta\otimes\hat R)((1\otimes\hat R(a_j))
      \alpha(x_j x_k^*)(1\otimes\hat R(a_k^*))\\
&=&\hat R(\eta\otimes\iota)(zz^*)\geq 0,
\end{eqnarray*}
where $z=\sum_j (1\otimes\hat R(a_j))\alpha(x_j)$.

To see that $y$ is invertible, note that
$(\iota\otimes\hat\varepsilon)(y)=1$. For this apply
$\iota\otimes\hat\varepsilon$ to the identity
$(\eta\otimes\iota)\alpha (b)=(\eta\otimes\iota)((b\otimes
1)(\iota\otimes\hat S)(y))$, use the property
$(\iota\otimes\hat\varepsilon)\alpha =\iota$ from Proposition
\ref{1.3} and faithfulness of $\eta$. Now as $y$ is self-adjoint,
we may write the cocycle identity as
$(\iota\otimes\hat\D)(y)=(y\otimes 1)(\alpha\otimes\iota)(y)$.
Next apply the map $(\iota\otimes m)(\iota\otimes\iota\otimes\hat
S)$ on both sides. This gives
\begin{eqnarray*}
1\otimes 1
&=&(1\otimes\hat\varepsilon(\cdot)1)(y)
    =(\iota\otimes m)(\iota\otimes\iota\otimes\hat S)
     (\iota\otimes\hat\D)(y)\\
&=&(\iota\otimes m)(\iota\otimes\iota\otimes\hat S)
      ((y\otimes 1)(\alpha\otimes\iota)(y))
     =y(\iota\otimes m)(\alpha\otimes\iota)(\iota\otimes\hat S)(y).
\end{eqnarray*}
Since $y^*=y$, we conclude that $y$ is invertible with two-sided
inverse $(\iota\otimes m)(\alpha\otimes\iota)(\iota\otimes\hat
S)(y)$.
\ep

Consider now the right coaction $\hat\Phi \colon A\to
M(A\otimes\hat A)$ introduced above.

\begin{prop} \label{1.12}
For the Haar state $\varphi$ on $A$, the following properties
hold:
\enu{i}
$(\sigma^\varphi_t\otimes\hat\tau_{-t})\hat\Phi
=\hat\Phi\sigma^\varphi_t$;
\enu{ii}
the state $\varphi$ is quasi-invariant with Radon-Nikodym cocycle
$y=W(1\otimes\rho^{-2})W^*$.
\end{prop}

\bp
Recall that $\hat\tau_t(x)=\rho^{it}x\rho^{-it}$ and $\hat\Phi
(a)=W(a\otimes 1)W^*$. By $\sigma^\varphi_t (a)=f_{it}*a*f_{it}$
and (\ref{e1.1}), we get $(\sigma^\varphi_t\otimes\iota)(W)
=(1\otimes\rho^{it})W(1\otimes\rho^{it})$. This implies (i).

\smallskip

For $a\in A$, we have
$$
(\varphi\otimes\iota)\hat\Phi (a)
=(\varphi\otimes\iota)(W(a\otimes 1)W^*)
=(\varphi\otimes\iota)
  ((a\otimes 1)(\iota\otimes\hat S)((\iota\otimes\hat S^{-1})(W^*)
  (\sigma^\varphi_{-i}\otimes\hat S^{-1})(W))),
$$
so $\varphi$ is quasi-invariant with Radon-Nikodym cocycle
$(\iota\otimes\hat S^{-1})(W^*)(\sigma^\varphi_{-i}\otimes\hat
S^{-1})(W)$. Because $(\iota\otimes\hat S^{-1})(W^*)=W$, we thus
must show that $(\sigma^\varphi_{-i}\otimes\hat
S^{-1})(W)=(1\otimes\rho^{-2})W^*$. We have
\begin{eqnarray*}
(\sigma^\varphi_{-i}\otimes\hat S^{-1})(W)
&=&(\iota\otimes\hat S^{-1})((1\otimes\rho)W
        (1\otimes\rho))
     =(1\otimes\rho^{-1})(\iota\otimes\hat S^{-1})(W)
        (1\otimes\rho^{-1})\\
&=&(1\otimes\rho^{-2})(1\otimes\rho)(\iota\otimes\hat S^{-1})(W)
        (1\otimes\rho^{-1})
=(1\otimes\rho^{-2})(1\otimes\hat S^2)(\iota\otimes\hat S^{-1})(W)\\
&=&(1\otimes\rho^{-2})W^*
\end{eqnarray*}
as desired. Note that the cocycle identity, positivity and
invertibility for $y$ follow from Propositions~\ref{1.10}
and~\ref{1.11}, but can easily be checked directly.
\ep

\bigskip\bigskip

\section{Random Walks on Discrete Quantum Groups} \label{s2}

\subsection{Quantum Markov Chains} \label{s2.1}

Let $(\hat A,\hat\D)$ be a discrete quantum group. For any bounded
linear functional $\phi$ on $\hat A$, we define the convolution
operator $P_\phi \colon M(\hat A)\to M(\hat A)$ by $P_\phi
=(\phi\otimes\iota)\hat\D$. Here $M(\hat A)$ is the multiplier
algebra of $\hat A$, so $M(\hat A)$ is the C$^*$-algebraic direct
product of $B(H_s)$, $s\in I$. Note that $P_\phi (\hat
A)\subset\hat A$ and if $\eta$ is another bounded linear
functional on $\hat A$, then $P_\phi P_\eta =P_{\phi\eta}$. If
$\phi$ is a state, then $P_\phi$ is a unital completely positive
map, and we call it the Markov operator associated with $\phi$.

We shall usually consider $\phi$ in the space $\C$ of
$\Phi$-invariant linear functionals on $\hat A$ introduced in
Subsection \ref{s1.4}. The following proposition partially
justifies this. Let $Z(\hat A)$ denote the center of $\hat A$, so
$Z(\hat A)$ is the C$^*$-algebra generated by $I_s$, $s\in I$, so
$Z(\hat A)\cong c_0 (I)$.

\begin{prop}
For a bounded linear functional $\phi$ on $\hat A$, we have
$P_\phi (Z(\hat A))\subset Z(\hat A)$ if and only if $\phi\in\C$.
\end{prop}

\bp
By Lemma \ref{1.2}, we see that $P_\phi(I_0)\in Z(\hat A)$ if and
only if $\phi (n^{\bar s}_{ij})=0$ for $i\ne j$ and the number
$\lambda_{\bar s}=f_{-1}(u^s_{ii})\phi (n^{\bar s}_{ii})$ is
independent of $i$, for any $s\in I$. Since by Lemma
\ref{1.1}(ii), we have $f_{-1}(u^s_{ii})=f_1 (v^{\bar s}_{ii})$,
these conditions mean that $\phi|_{B(H_{\bar s})}=d_{\bar
s}\lambda_{\bar s}\phi_{\bar s}$, so $\phi\in\C$. We have thus
proved the forward implication. To prove that $P_\phi (Z(\hat
A))\subset Z(\hat A)$ whenever $\phi\in\C$, note that by the above
argument, we see that $P_{\phi_s}(I_0 )=\frac{1}{d_s^2}I_{\bar s}$
and therefore $P_{\phi}(I_0)\in Z(\hat A)$ for all $\phi\in\C$.
Also $P_\phi (I_s )=d_s^2 P_\phi P_{\phi_{\bar s}}(I_0)=d_s^2
P_{\phi\phi_{\bar s}}(I_0)$, and as $\C$ is an algebra by
Lemma~\ref{1.8}, we see that $P_\phi (I_s )\in Z(\hat A)$.
\ep

Suppose $\phi\in\C$ is a state. As in \cite{B1} we can construct
the associated quantum Markov chain (with initial distribution
$\hat\varepsilon$ and generator $P_\phi$) in the sense of
\cite{AFL}. It consists of:
\enu{i}
a von Neumann algebra $M(\hat A)^\infty$ and a normal state
$\phi^\infty$ given by $\otimes_{-\infty}^{-1}(M(\hat A),\phi)$;
\enu{ii}
unital $*$-homomorphisms $j_k \colon M(\hat A)\to M(\hat
A)^\infty$ given by $j_k(a) =\dots\otimes1\otimes\hat\D^{k-1}(a)$
for $k\geq 1$ and $a\in M(\hat A)$, and $j_0 =\hat\varepsilon$.

Here $\hat\D^k$ is defined inductively by $\hat\D^0 =\iota$,
$\hat\D^1=\hat\D$ and $\hat\D^{k+1}=(\hat\D\otimes\iota)\hat\D^k$.
Using $\phi=\hat\varepsilon P_\phi$, the crucial property
\begin{equation} \label{e2.1}
\phi^\infty (j_0 (a_0)\dots j_k(a_k))
=\hat\varepsilon (a_0 P_\phi(\dots P_\phi(a_{k-1}P_\phi
(a_k))\dots ))
\end{equation}
is then easily checked.

The algebra $M(\hat A)^\infty$ should be thought of as the algebra
of measurable functions on the path space of our quantum Markov
chain. Since $P_\phi (Z(\hat A))\subset Z(\hat A)=c_0 (I)$, the
operator $P_\phi$ determines a classical Markov chain on $I$ with
kernel $\{p_\phi(s,t)\}_{s,t\in I}$, so
$P_\phi(I_t)I_s=p_\phi(s,t)I_s$. Let $(\Omega ,\7P_0 )$ be the
corresponding path space. Thus
$\displaystyle\Omega=\prod^{-1}_{-\infty}I$ and the measure
$\7P_0$ is defined on cylindrical sets as follows:
\begin{equation} \label{e2.2}
\7P_0(\{\omega\in\Omega\,|\,\omega_{-n}
  =s_{-n},\dots,\omega_{-1}=s_{-1}\})
=p_\phi(0,s_{-1})p_\phi(s_{-1},s_{-2})\dots
p_\phi(s_{-n+1},s_{-n}).
\end{equation}

The following proposition is essentially from \cite{B1}.

\begin{prop}
There exists an embedding $j^\infty\colon
L^\infty(\Omega,\7P_0)\hookrightarrow(M(\hat
A)^\infty,\phi^\infty)$ uniquely determined by
$j^\infty(a_{-n}\otimes\dots\otimes a_{-1})=j_n(a_{-n})\dots
j_1(a_{-1})$ for $a_{-n},\dots,a_{-1}\in Z(\hat A)$.
\end{prop}

\bp
We first prove that $j_{k+1}(a)$ commutes with $j_{l+1}(b)$ for
arbitrary $k,l\ge0$ and $a,b\in Z(\hat A)$. We can assume that
$l=k+n$ for some $n\in\7N$. Thus we need to show that
$1\otimes\dots\otimes 1\otimes\hat\D^k(a)$ and $\hat\D^{k+n}(b)$
commute in $M(\otimes^{-1}_{-k-n}\hat A)$. This is true because
$1\otimes\dots\otimes
1\otimes\hat\D^k(a)=(\iota\otimes\dots\otimes
\iota\otimes\hat\D^k) (1\otimes\dots\otimes 1\otimes a)$,
$\hat\D^{k+n}(b) =(\iota\otimes\dots\otimes
\iota\otimes\hat\D^k)\hat\D^n(b)$ by coassociativity of $\hat\D$,
and because $1\otimes\dots\otimes 1\otimes a$ commutes with
$\hat\D^n(b)$. So there exists a unital $*$-homomorphism
$j^\infty\colon\otimes^{-1}_{-\infty}l^\infty(I)\to M(\hat
A)^\infty$ determined by $j^\infty(a_{-n}\otimes\dots\otimes
a_{-1})=j_n(a_{-n})\dots j_1(a_{-1})$ for $a_{-n},\dots,a_{-1}\in
Z(\hat A)$. Using equalities~(\ref{e2.1}) and~(\ref{e2.2}) it is
easy to check that $\phi^\infty j^\infty=\7P_0$. Thus $j^\infty$
extends to the required normal embedding.
\ep

\subsection{Transience}

\begin{defin} Suppose $\phi\in\C$ is a state.
\enu{i}
We say that $\phi$ is transient if the corresponding classical
random walk on $I$ is transient, that is, if the sum
$\sum_{n=0}^\infty p_{\phi^n}(s,t)$ is finite for all $s,t\in I$.
\enu{ii}
By $\supp\phi$ we mean the set $\{s\in I\ |\ \phi (I_s )\ne 0\}$.
\enu{iii}
We say that $\phi$ is generating if for any $s\in I$, there exists
$n\in\7N$ such that $\phi^n (I_s)>0$, that is $\cup_n\supp\phi^n
=I$.
\end{defin}

Define an anti-linear isometric operator $\phi\mapsto\check\phi$
on $\C$ by $\check\phi_s =\phi_{\bar s}$ for $s\in I$.

Let $U$ be a finite dimensional unitary corepresentation of
$(A,\D)$. Consider the state $\phi_U$ given by
$\phi_U=d_U^{-1}\Tr\,\pi_U(\cdot\rho^{-1})$, where
$d_U=\Tr\,\pi_U(\rho^{-1})=\Tr\,\pi_U(\rho)$
is the quantum dimension. Then $\phi_U\in\C$. Define $N_{U,s}^t$
to be the multiplicity of $U^t$ in $U\times U^s$.

\begin{lem} \label{2.4}
With the above notation the following properties hold:
\enu{i}
$\phi_U\phi_V=\phi_{U\times V}$ for any finite dimensional unitary
corepresentations $U$ and $V$;
\enu{ii}
$\check\phi_U=\phi_{\bar U}$;
\enu{iii}
the mapping $\phi\mapsto\check\phi$ on $\C$ is anti-multiplicative;
\enu{iv}
$p_{\phi_U}(s,t)=\frac{d_t}{d_U d_s}N^t_{U,s}$;
\enu{v}
$p_{\check \phi}(s,t)=(\frac{d_t}{d_s})^2 p_\phi (t,s)$ for any
state $\phi\in\C$.
\end{lem}

\bp
Parts (i) and (ii) follow immediately from definitions. Part (iii)
follows from (i) and (ii) and the property $\overline{U\times
V}\cong\bar V\times\bar U$. For part (iv) apply $\phi_s$ to the
identity $P_{\phi_U}(I_t)I_s=p_{\phi_U}(s,t)I_s$. This yields (see
the proof of Lemma \ref{1.8})
$$
p_{\phi_U}(s,t)=\phi_s P_{\phi_U}(I_t)=(\phi_U\phi_s )(I_t)
=\phi_{U\times U^s}(I_t)=\sum_w \frac{d_w}{d_U d_s}N_{U,s}^w
\phi_w (I_t) =\frac{d_t}{d_U d_s}N_{U,s}^t .
$$
Part (v) now follows for $\phi=\phi_U$ by the Frobenius
reciprocity $N^t_{U,s}=N^s_{\bar U,t}$, and for general $\phi$ by
linearity.
\ep

\begin{cor} \label{2.5}
A state $\phi\in\C$ is generating if and only if the classical
random walk is irreducible in the sense that for all $s,t\in I$,
there exists $n\in\7N$ such that $p_{\phi^n}(s,t)>0$.
\end{cor}

\bp
Since $p_\phi (0,s)=\phi (I_s)$, we see that irreducibility of the
classical random walk implies the generating property for $\phi$.
Conversely, suppose $\phi$ is generating. Since
$(\check\phi)^n=(\phi^n)\,\check{}$ by Lemma~\ref{2.4}(iii), we
have $\supp(\check\phi)^n=\overline{\supp\phi^n}$, so $\check\phi$
is generating. Now given $s,t\in I$ we can find $k,l\in\7N$ such
that $(\check\phi)^k(I_s)>0$ and $\phi^l(I_t)>0$. By
Lemma~\ref{2.4}(v) we thus get
$p_{\phi^k}(s,0)=d_s^{-2}p_{(\check\phi)^k}(0,s)>0$ and
$p_{\phi^l}(0,t)>0$. Since
$$
p_{\phi^{k+l}}(s,t)=\sum_{w\in I}p_{\phi^k}(s,w)p_{\phi^l} (w,t),
$$
we conclude that $p_{\phi^{k+l}}(s,t)>0$.
\ep

Note also that Lemma \ref{2.4} implies that $\phi$ is generating
if and only if for any $s\in I$ there exist
$s_1,\dots,s_n\in\supp\phi$ such that $U^s$ is a
subcorepresentation of $U^{s_1}\times\dots\times U^{s_n}$. In
particular, if $(A,\D)$ is a compact matrix pseudogroup with
fundamental corepresentation $U$ \cite{woro1}, the state $\phi_U$
is generating.

\medskip

In the classical theory it is usually difficult to check the
transience condition. The following result shows, however, that in
the generic quantum group case transience is automatic.

\begin{thm} \label{2.6}
Suppose $\phi\in\C$ is a state for which there exists
$w\in\supp\phi$ with $\pi_w(\rho)\ne I_w$ (or equivalently, $\dim
H_w<d_w$). Then there exists $\lambda<1$ such that
$p_{\phi^n}(s,t)\le\frac{d_t\dim H_s}{d_s\dim H_t}\lambda^n$ for
any $s,t\in I$ and $n\in\7N$. In particular, $\phi$ is transient.
\end{thm}

\bp
Recall that since $\Tr\,\pi_s(\rho)=\Tr\,\pi_s(\rho^{-1})$, we
have by Schwarz inequality
$$
\dim H_s=\Tr\,\pi_s(\rho^\2\rho^{-\2})\le d_s,
$$
with equality if and only if $\pi_s(\rho^\2)=\pi_s(\rho^{-\2})$,
that is $\pi_s(\rho)=I_s$.

Let $\phi=\sum_r\lambda_r\phi_r$, where $\lambda_r\ge0$ and
$\sum_r\lambda_r=1$. By Lemma \ref{2.4}(i) we have
$$
\phi^n=\sum_{r_1,\dots,r_n}\lambda_{r_1}\dots\lambda_{r_n}
\phi_{U^{r_1}\times\dots\times U^{r_n}}.
$$
Thus Lemma \ref{2.4}(iv) and the inequality
$N_{U,s}^t\le\frac{\dim H_U\dim H_s}{\dim H_t}$ entails
\begin{eqnarray*}
p_{\phi^n}(s,t)
&=&\sum_{r_1,\dots,r_n}\lambda_{r_1}\dots\lambda_{r_n}
       \frac{d_t}{d_{r_1}\dots d_{r_n}d_s}
      N_{U^{r_1}\times\dots\times U^{r_n},s}^t\\
&\le&\sum_{r_1,\dots,r_n}\lambda_{r_1}\dots\lambda_{r_n}
       \frac{\dim H_{r_1}}{d_{r_1}}\dots
          \frac{\dim H_{r_n}}{d_{r_n}}
       \frac{d_t\dim H_s}{d_s\dim H_t}
  =\lambda^n\frac{d_t\dim H_s}{d_s\dim H_t},
\end{eqnarray*}
where $\lambda=\sum_r\lambda_r\frac{\dim H_r}{d_r}<1$.
\ep

In particular, if $(A,\D)$ is a $q$-deformation of (the algebra of
continuous functions on) a semisimple compact Lie group (see e.g.
\cite{LS}), then any state $\phi\ne\hat\varepsilon$ in $\C$ is
transient. The assumptions of Theorem~\ref{2.6} are also satisfied
for any non-Kac algebra $(\hat A,\hat\D)$ with a generating state
$\phi\in\C$.

\medskip

Our next goal is to extend $P_\phi$ to a larger subspace of
$\A'=M(\hat\A)$. We consider $\A'$ with weak$^*$ topology,
which coincides with the Tikhonov topology on
the algebraic direct product
$M(\hat\A)=\prod_s B(H_s)$. So $M(\hat\A)$ is a
complete locally convex space. Let $\I$ be the collection of all
finite subsets of $I$. For $X\subset I$ denote by $F_X$ the
projection from $M(\hat\A)$ onto $\prod_{s\in X}B(H_s)$. So a net
$\{x_i\}_i$ converges to $x$ in $M(\hat\A)$ if and only if
$F_X(x_i)\to F_X(x)$ for all $X\in\I$.

Let $\phi$ be a positive linear functional in $\C$. Denote by
$D(P_\phi)_+$ the set of all $x\in M(\hat\A)_+$ such that the net
$\{P_\phi F_X(x)\}_{X\in\I}$ in $M(\hat\A)$ is Cauchy, or
equivalently, such that the set $\{\pi_sP_\phi F_X(x)\}_{X\in\I}$
is bounded for any $s\in I$. For $x\in D(P_\phi)_+$ we set
$P_\phi(x) =\lim_{X\in\I} P_\phi F_X(x)$, and then form the linear
extension $P_\phi \colon D(P_\phi)\to M(\hat\A)$, where
$D(P_\phi)$ is the linear span of $D(P_\phi)_+$. Note that
$D(P_\phi)\cap M(\hat\A)_+ =D(P_\phi)_+$. If the support of $\phi$
is finite, then $D(P_\phi)=M(\hat\A)$. Moreover, as $\phi\in\A$ in
this case, the formula $P_\phi =(\phi\otimes\iota)\hat\D$ is
meaningful.

\begin{lem}\label{2.7}
\mbox{}
\enu{i}
If $0\leq x\leq y$ for $x\in M(\hat\A)$ and $y\in D(P_\phi)$, then
$x\in D(P_\phi)$ and $0\leq P_\phi(x)\leq P_\phi(y)$.
\enu{ii}
If $x_i\to x\in M(\hat\A)$ and $0\leq x_i\leq y$ for some $y\in
D(P_\phi)$, then $x\in D(P_\phi)$ and $P_\phi(x_i)\to P_\phi(x)$.
\end{lem}

\bp
Assertion (i) is obvious from definitions. As for (ii), note that
$x\in D(P_\phi )$ by (i). Fix $X\in\I$ and $\varepsilon>0$. Then
there exists $Y\in\I$ such that $\|F_X P_\phi (y-F_Y
(y))\|<\varepsilon$.  There also exists $i_0$ such that $\|F_X
P_\phi F_Y(x-x_i)\|<\varepsilon$ for all $i\geq i_0$. Thus
$$
\|F_X(P_\phi (x)-P_\phi({x_i}))\|
\leq\|F_X P_\phi (x-F_Y (x))\|+\|F_X P_\phi (x_i-F_Y (x_i))\|
   +\|F_X P_\phi F_Y (x-x_i))\|<3\varepsilon,
$$
since $0\leq F_XP_\phi (x-F_Y (x))\leq F_X P_\phi (y-F_Y (y))$ and
$0\leq F_XP_\phi (x_i-F_Y (x_i))\leq F_X P_\phi (y-F_Y (y))$.
\ep

We proceed to define the potential operator $G_\phi
=\sum_{n=0}^\infty P_\phi^n\colon D(G_\phi)\to M(\hat\A)$. By
definition its domain $D(G_\phi)\subset M(\hat\A)$ is the linear
span of $D(G_\phi)_+$, where $x\in D(G_\phi)_+$ whenever $x\geq
0$, $x\in \cap_n D(P^n_\phi)$ and the series $\sum_{n=0}^\infty
P_\phi^n (x)$ converges.

\medskip

By definition $\phi$ is transient if and only if $I_s\in
D(G_\phi)$ for all $s\in I$, if and only if $\hat\A\subset
D(G_\phi )$. Thus in this case the domain of $G_\phi$ is dense,
but $G_\phi$ is by no means continuous.

Note also that if $\phi$ is transient, then $G_\phi(\hat\A)\subset
M(\hat A)$ by complete maximum principle (see 
e.g.~\cite[Theorem~2.1.12]{R}).

\begin{defin}
We say that an element $x\in M(\hat\A)$ is
\enu{i}
harmonic if $P_\phi (x)=x$;
\enu{ii}
superharmonic if $x\ge0$ and $P_\phi (x)\leq x$;
\enu{iii}
a potential of an element $y\in D(G_\phi)_+$ if $x=G_\phi (y)$.
\end{defin}

Recall that there exists a strong connection between transience
and superharmonic elements: $\phi$~is transient if and only if
there exist non-constant central (that is, lying in $Z(\A')$)
superharmonic elements (see e.g. \cite[Theorem~1.16]{W}).

\begin{lem} \label{2.9}
\mbox{\ }
\enu{i}
If $x\in D(G_\phi)$, then $G_\phi (x)\in D(P_\phi)$ and $P_\phi
G_\phi (x)=G_\phi(x)-x$. In particular, any potential is
superharmonic.
\enu{ii} A superharmonic element $x$ is a potential if and only if
$P^n_\phi(x)\to0$ in $M(\hat\A)$. In particular, any superharmonic
element which is majorized by a potential, is a potential itself.
\end{lem}

\bp
In proving (i) we can suppose that $x\ge0$. For any $Y\in\I$ the
set $\{F_YP_\phi F_XG_\phi(x)\}_X$ is bounded. Indeed, if
$x_n=\sum_{k=0}^n P^k_\phi (x)$, then $F_YP_\phi F_Xx_n\nearrow
F_YP_\phi F_XG_\phi(x)$ and $F_YP_\phi F_Xx_n\le
F_YP_\phi(x_n)=F_Y(x_{n+1}-x)\le F_YG_\phi(x)$. Hence
$G_\phi(x)\in D(P_\phi)$. Then by Lemma \ref{2.7}(ii) we deduce
$\displaystyle P_\phi G_\phi(x)=\lim_{n\to\infty}P_\phi (x_n)
=\lim_{n\to\infty}(x_{n+1}-x)=G_\phi (x)-x$.

To prove (ii) consider a potential $x=G_\phi(y)$. Then by (i)
$P_\phi^n(x)=G_\phi(y)-y_{n-1}$, where $y_{n-1}=\sum_{k=0}^{n-1}
P^k_\phi (y)$. Thus $P^n_\phi(x)\to0$. Conversely, if
$P_\phi(x)\le x$ and $P^n_\phi(x)\to0$, then $x=G_\phi(y)$ with
$y=x-P_\phi(x)$, since $\sum_{k=0}^n P^k_\phi
(y)=x-P^{n+1}_\phi(x)\to x$.
\ep

\subsection{Balayage Theorem}

The aim of this subsection is to show that any superharmonic
element can be  canonically approximated by potentials. The proof
is essentially the same as in the classical theory. However, since
most classical proofs use some simplifications arising from
probabilistic interpretations and the fact the Markov operator
acts on functions (so the space of self-adjoint elements is a
lattice), we present a detailed argument.

For a positive linear functional $\phi\in\C$ and $Y\in\I$,
consider the linear operator
$$
P_{\phi, Y}=\sum^\infty_{n=0}[(\iota- F_Y)P_\phi]^nF_Y.
$$
Note that if $F_Y(x)\in D(G_\phi)_+$, then $x\in D(P_{\phi, Y})$
and $P_{\phi, Y}(x)\le G_\phi F_Y(x)$.

Define also
$$
G^Y_\phi=\sum^\infty_{n=0}[(\iota- F_Y)P_\phi]^n(\iota-F_Y).
$$
Then $D(G_\phi)\subset D(G^Y_\phi)$ and $G^Y_\phi(x)\le G_\phi(x)$
for $x\in D(G_\phi)_+$.

\begin{lem} \label{2.10}
For any $x\in D(G_\phi)$ and $Y\in\I$, we have $G_\phi(x)\in
D(P_{\phi, Y})$ and
$$
G_\phi(x)=G^Y_\phi(x)+P_{\phi, Y}G_\phi(x).
$$
\end{lem}

\bp
Let $x\ge0$. Set $x_n=\sum^n_{k=0}P^k_\phi(x)$. Note that
$$
P^k_\phi=\sum^k_{m=0}[(\iota- F_Y)P_\phi]^{k-m}F_Y P^m_\phi
+[(\iota- F_Y)P_\phi]^k(\iota-F_Y).
$$
This can be verified by induction as follows:
\begin{eqnarray*}
(\iota-F_Y)P^k_\phi
&=&(\iota-F_Y)P_\phi(\sum^{k-1}_{m=0}
       [(\iota- F_Y)P_\phi]^{k-m-1}F_Y P^m_\phi
    +[(\iota- F_Y)P_\phi]^{k-1}(\iota-F_Y))\\
&=&\sum^{k-1}_{m=0}[(\iota- F_Y)P_\phi]^{k-m}F_Y P^m_\phi
+[(\iota- F_Y)P_\phi]^k(\iota-F_Y).
\end{eqnarray*}
Consequently
\begin{eqnarray*}
x_n&=&\sum^n_{k=0}\sum^k_{m=0}
        [(\iota- F_Y)P_\phi]^{k-m}F_Y P^m_\phi(x)
+\sum^n_{k=0}[(\iota- F_Y)P_\phi]^k(\iota-F_Y)(x)\\
&=&\sum^n_{k=0}[(\iota- F_Y)P_\phi]^kF_Y\sum^{n-k}_{m=0}P^m_\phi(x)
+\sum^n_{k=0}[(\iota- F_Y)P_\phi]^k(\iota-F_Y)(x)\\
&=&\sum^n_{k=0}[(\iota- F_Y)P_\phi]^kF_Y(x_{n-k})
+\sum^n_{k=0}[(\iota- F_Y)P_\phi]^k(\iota-F_Y)(x).
\end{eqnarray*}
This identity already shows that $G_\phi(x)\in D(P_{\phi, Y})$.
Indeed, since $G_\phi(x)\in D(P_\phi^k)$ for any $k\in\7N$ by
Lemma \ref{2.9}, it follows by Lemma \ref{2.7}(ii) that
$[(\iota-F_Y)P_\phi]^kF_Y(x_{n-k})
\nearrow[(\iota-F_Y)P_\phi]^kF_YG_\phi(x)$
when $n\to\infty$. Thus the identity above implies
$\sum^m_{k=0}(\iota-F_Y)P_\phi]^kF_YG_\phi(x)\le G_\phi(x)$ for
any $m\in\7N$, whence $G_\phi(x)\in D(P_{\phi, Y})$.

To complete the proof of Lemma it is enough to show that
$\sum^n_{k=0}[(\iota- F_Y)P_\phi]^kF_Y(x_{n-k}-G_\phi(x))\to0$ as
$n\to\infty$. Since $G_\phi(x)-x_{n-k}\le G_\phi(x)$ and
$G_\phi(x)\in D(P_{\phi,Y})$, there exists $n_0\in\7N$ such that
$\sum^n_{k=n_0}[(\iota- F_Y)P_\phi]^kF_Y(x_{n-k}-G_\phi(x))$ is
close to zero for all $n\ge n_0$. On the other hand, as we have
already remarked, we have
$[(\iota-F_Y)P_\phi]^kF_Y(x_{n-k}-G_\phi(x))\to0$ for fixed $k$ as
$n\to\infty$. This ends the proof.
\ep

\begin{thm} \label{2.11}
Suppose $\phi\in\C$ is transient and consider a superharmonic
element $x\in M(\hat\A)$. Then for any $Y\in\I$, the following
properties hold:
\enu{i}
$P_{\phi, Y}(x)\le x$ and $F_YP_{\phi, Y}(x)=F_Y(x)$;
\enu{ii}
if $y$ is a superharmonic element satisfying $F_Y(x)\le F_Y(y)$,
then $P_{\phi, Y}(x)\le y$;
\enu{iii}
$P_{\phi, Y}(x)$ is a potential.
\end{thm}

Thus any superharmonic element $x$ is approximated by potentials
$P_{\phi, Y}(x)$ in the topology of~$M(\hat\A)$. Moreover,
$P_{\phi, Y}(x)$ is the smallest element in the set of
superharmonic elements majorizing $x$ on $Y$.

\smallskip

\bpp{Theorem \ref{2.11}}
If $x$ is a potential, say $x=G_\phi(y)$, then
$$
x=G_\phi(y)=G^Y_\phi(y)+P_{\phi, Y}G_\phi(y)
\ge P_{\phi, Y}G_\phi(y)=P_{\phi, Y}(x).
$$
For general superharmonic $x$, consider the operator
$P_{\lambda\phi}=\lambda P_\phi$ for $0<\lambda<1$. Clearly $x$ is
also superharmonic with respect to $P_{\lambda\phi}$. Since
$P_{\lambda\phi}^n(x)\le\lambda^n x\to0$, by Lemma 2.9(ii) we
deduce that~$x$ is a potential with respect to $P_{\lambda\phi}$.
Hence $P_{\lambda\phi,Y}(x)\le x$. This is equivalent to the fact
that $\sum^n_{k=0}\lambda^k[(\iota- F_Y)P_\phi]^kF_Y(x)\le x$ for
all $n\in\7N$ and $\lambda<1$. Letting $\lambda\to1-0$ and
$n\to\infty$, we conclude that $P_{\phi,Y}(x)\le x$. The equality
$F_YP_{\phi,Y}(x)=F_Y(x)$ follows by definition. Thus (i) is
proved.

\smallskip

If $y\ge0$ is such that $P_\phi(y)\le y$ and $F_Y(x)\le
F_Y(y)$, we get
$$
P_{\phi,Y}(x)=P_{\phi,Y}F_Y(x)\le P_{\phi,Y}F_Y(y)
=P_{\phi,Y}(y)\le y,
$$
where the last inequality follows from (i) applied to $y$. This
proves (ii).

\smallskip

To prove (iii), we shall first check that $P_{\phi,Y}(x)$ is
superharmonic. Since $P_{\phi,Y}(x)\le x$, we get
$$
F_YP_\phi P_{\phi,Y}(x)\le F_YP_\phi(x)\le F_Y(x)=F_YP_{\phi,Y}(x).
$$
The same proof as in Lemma \ref{2.9}(i) (applied to
$(\iota-F_Y)P_\phi$ and element $F_Y(x)$ instead of $P_\phi$ and
$x$) shows that $(\iota-F_Y)P_\phi
P_{\phi,Y}(x)=P_{\phi,Y}(x)-F_Y(x)$, so $(\iota-F_Y)P_\phi
P_{\phi,Y}(x)=(\iota-F_Y)P_{\phi,Y}(x)$. Thus
$$
P_\phi P_{\phi,Y}(x)=F_YP_\phi P_{\phi,Y}(x)
  +(\iota-F_Y)P_\phi P_{\phi,Y}(x)
\le F_YP_{\phi,Y}(x)+(\iota-F_Y)P_{\phi,Y}(x)=P_{\phi,Y}(x),
$$
so $P_{\phi,Y}(x)$ is superharmonic. Hence to
prove that it is a potential, by Lemma~\ref{2.9}(ii) it is enough
to show that it is majorized by a potential. But we obviously have
$P_{\phi,Y}(x)\le G_\phi F_Y(x)$. Thus (iii)
is also proved.
\ep

The previous proof is quite formal. It is applicable to any
positive operator on a complete ordered vector space with a given
increasing net $\I$ of ordered subspaces. We leave it to the
interested reader to formulate the precise setting for such a
generalization. It is worth noting that the assumption on the
subspaces in $\I$ to be monotonically complete, simplifies the
proof but is not at all necessary.

\subsection{Poisson Boundary and 0-2 Law}

Let $\phi\in\C$ be a generating state. Following Izumi \cite{I} we
denote by $H^\infty(M(\hat A),P_\phi)$ the set of bounded (that
is, the elements belonging to $M(\hat A)$) harmonic elements
with respect to $P_\phi$ and call it the Poisson boundary of
$(M(\hat A),P_\phi)$. The operator system $H^\infty(M(\hat
A),P_\phi)\subset M(\hat A)$ has a unique structure of a
C$^*$-algebra. To distinguish the product on $H^\infty(M(\hat
A),P_\phi)$ from the one on~$M(\hat A)$, we shall write $x\cdot y$
for the product of two harmonic elements $x$ and $y$. It is proved
in \cite{I} that
$$
x\cdot y=\lim_{n\to\infty}P_\phi^n(xy),
$$
where the limit is in the topology of $M(\hat\A)$. In fact,
$M(\hat A)$ is a von Neumann algebra, and $H^\infty(M(\hat
A),P_\phi)$ is a weakly operator closed subspace. Thus
$H^\infty(M(\hat A),P_\phi)$ is a von Neumann algebra.

\smallskip

Following \cite{I} and \cite{K} we shall presently give several
other descriptions of $H^\infty(M(\hat A),P_\phi)$.

\smallskip

We say that a sequence $\{x_n\}^\infty_{n=1}\subset M(\hat A)$ is
harmonic if $F_nP_\phi(x_{n+1})=x_n$, where $F_n=F_{\supp\phi^n}$
is the projection onto the C$^*$-algebra product of $B(H_s)$,
$s\in\supp\phi^n$. Note that the homomorphisms~$j_n$ introduced in
Subsection \ref{s2.1} fulfill $\phi^\infty j_n=\phi^n$ and
\begin{equation} \label{e2.3}
E_nj_{n+1}=j_n P_\phi,
\end{equation}
where $E_n\colon(M(\hat
A)^\infty,\phi^\infty)\to\otimes^{-1}_{-n}(M(\hat A),\phi)$ is the
$\phi^\infty$-preserving conditional expectation. Note also that
since $\phi^\infty j_n=\phi^n$, we have $j_n=j_nF_n$ and
$j_n|_{F_n(M(\hat A))}$ is injective. Thus if
$\{x_n\}^\infty_{n=1}$ is a bounded harmonic sequence, the
sequence $\{j_n(x_n)\}^\infty_{n=1}$ is a martingale, so it
converges in strong$^*$ operator topology to an element $x\in
M(\hat A)^\infty$ such that $E_n(x)=j_n(x_n)$ for all $n\in\7N$.
Conversely, if $x\in M(\hat A)^\infty$ satisfies
$E_n(x)\in\Im\,j_n$, and $x_n\in F_n(M(\hat A))$ is the element
uniquely determined by $E_n(x)=j_n(x_n)$, then
$\{x_n\}^\infty_{n=1}$ is a bounded harmonic sequence.

Any bounded harmonic element $x\in M(\hat A)$ defines a bounded
harmonic sequence $x_n=F_n(x)$, as $F_nP_\phi=F_nP_\phi F_{n+1}$,
which follows from the equality $\phi^nP_\phi=\phi^{n+1}$.
Conversely, if  $\{x_n\}^\infty_{n=1}$ is a bounded harmonic
sequence with $F_n(x_m)=F_m(x_n)$ for all $n,m\in\7N$, then the
unique element $x\in M(\hat A)$ such that $F_n(x)=x_n$ for all
$n\in\7N$ is harmonic, because
$$
F_nP_\phi(x)=F_nP_\phi F_{n+1}(x)
=F_nP_\phi(x_{n+1})=x_n=F_n(x)\ \ \forall n\in\7N,
$$
and $\cup_n\supp\phi^n=I$. It turns out that the assumption
$F_n(x_m)=F_m(x_n)$ is automatically fulfilled by the $0$-$2$ law.
This law was first proved by Ornstein and Sucheston \cite{OS}. The
proof was later clarified by Foguel \cite{F}. The same line of
arguments works also in the context of non-commutative
probability.

\begin{prop} \label{2.12}
Consider a unital positive map $P\colon A\to A$ on a C$^*$-algebra
$A$. Suppose there exist $m,k\in\7N$ and a positive map $S\colon
A\to A$ such that $S(1)$ is invertible, $P^{m+k}\ge S$ and $P^m\ge
S$. Then $\displaystyle\lim_{n\to\infty}\|P^{n+k}-P^n\|=0$.
\end{prop}

The name '$0$-$2$ law' is due to the fact that if $A$ is an
abelian von Neumann algebra, then the existence of $S$ means
precisely that $\|P^{m+k}-P^m\|<2$. Hence, if $A$ is an abelian
C$^*$-algebra, then either $\|P^{n+k}-P^n\|=2$ for all $n\in\7N$,
or $\displaystyle\lim_{n\to\infty}\|P^{n+k}-P^n\|=0$.

\smallskip

\bpp{Proposition \ref{2.12}}
Set $h=m+k$. We claim that there exist positive maps $S_{ij}$ and
$T_j$ on~$A$ such that $S_{ij}(1)$ is invertible and
\begin{equation} \label{e2.4}
P^{ijh}=S_{ij}(\iota+P^k)^j+T^i_j\ \ \hbox{for all}\ i,j\in\7N.
\end{equation}
To this end, take $S_{11}=\2 S$ and $T_1=P^h-\2 S(\iota+P^k)=\2
(P^h-S)+\2 (P^m-S)P^k$. We define $S_{1j}$ and $T_j$ by induction
on $j$ using the equality
$$
P^{(j+1)h}=P^{jh}P^h=S_{1j}P^h(\iota+P^k)^j+T_jP^h
=S_{1j}S_{11}(\iota+P^k)^{j+1}+(S_{1j}T_1(\iota+P^k)^j+T_jP^h),
$$
so $S_{1,j+1}=S_{1j}S_{11}$ and
$T_{j+1}=S_{1j}T_1(\iota+P^k)^j+T_jP^h$. Then we define $S_{ij}$
by induction on $i$ using
$$
P^{(i+1)jh}=S_{ij}P^{jh}(\iota+P^k)^j+T_j^iP^{jh}
=(S_{ij}P^{jh}+T_j^iS_{1j})(\iota+P^k)^j+T_j^{i+1},
$$
which proves our claim.

Applying (\ref{e2.4}) to the unit, we conclude that $\|T_j\|<1$
and $\|S_{ij}\|\le 2^{-j}$. Equation (\ref{e2.4}) also yields
$$
\|P^{ijh}(\iota-P^k)\|\le\|S_{ij}(\iota+P^k)^j(\iota-P^k)\|
 +\|T_j^i(\iota-P^k)\|.
$$
The first term on the right hand side converges to zero as
$j\to\infty$ uniformly in $i\in\7N$, since
\begin{eqnarray*}
\|S_{ij}(\iota+P^k)^j(\iota-P^k)\|
&=&\|S_{ij}(\iota+\sum^j_{r=1}({j\choose r}-{j\choose r-1})
     P^{rk}-P^{(j+1)k})\|\\
&\le&2^{-j}(2+\sum^j_{r=1}|{j\choose r}-{j\choose r-1}|),
\end{eqnarray*}
and the latter expression converges to zero by known properties of
binomial coefficients \cite{OS}. On the other hand, for fixed $j$,
the second term converges to zero as $i\to\infty$, since
$\|T_j^i(\iota-P^k)\|\le2\|T_j\|^i$. We see that under an
appropriate choice of $i$ and $j$ the norm
$\|P^{ijh}(\iota-P^k)\|$ can be made arbitrarily small. Since the
sequence $\{\|P^n(\iota-P^k)\|\}_n$ is decreasing, its limit is
therefore zero.
\ep

Let us return to the proof of $F_n(x_m)=F_m(x_n)$. If
$\supp\phi^n\cap\supp\phi^m=\emptyset$, then
$F_n(x_m)=F_m(x_n)=0$. If $s\in\supp\phi^n\cap\supp\phi^m$, then
$P_\phi^m\ge cP_{\phi_s}$ and $P_\phi^n\ge cP_{\phi_s}$ for $c>0$
such that $c\phi_s\le\phi^m$ and $c\phi_s\le\phi^n$. Thus if
$l=m-n\ge0$, then $\|P^{k+l}_\phi-P^k_\phi\|\to0$ as $k\to\infty$.
Since $F_n(x_m)=F_nF_mP_\phi^k(x_{m+k})$ and
$F_m(x_n)=F_mF_nP_\phi^{k+l}(x_{n+l+k})=F_nF_mP_\phi^{k+l}(x_{m+k})$,
we see that
$\|F_n(x_m)-F_m(x_n)\|\le\|P^{k+l}_\phi-P^k_\phi\|\,\|x_{m+k}\|\to0$
as $k\to\infty$.

Summarizing the discussion above we obtain (see \cite{I,K})

\begin{thm} \label{2.13}
The following linear spaces are canonically isomorphic:
\enu{i}
the space $H^\infty(M(\hat A),P_\phi)$ of bounded harmonic
elements;
\enu{ii}
the space of bounded harmonic sequences;
\enu{iii}
the space of elements $x\in M(\hat A)^\infty$ such that
$E_n(x)\in\Im\,j_n$ for all $n\in\7N$.

Explicitly, the correspondence $\theta$ between (i) and (iii)
associates to each $x\in H^\infty(M(\hat A),P_\phi)$ an element
$\theta(x)\in M(\hat A)^\infty$ uniquely determined by
$E_n\theta(x)=j_n(x)$.
\epp
\end{thm}

Since $E_n(\Im\,j_{n+1})\subset\Im\,j_n$, the space described in
part (iii) is, in fact, a von Neumann subalgebra of $M(\hat
A)^\infty$. As $\theta$ is a unital completely positive and
isometric map of $H^\infty(M(\hat A),P_\phi)$ onto this
subalgebra, it is a $*$-homomorphism. Thus $\theta$ is an
embedding of $H^\infty(M(\hat A),P_\phi)$ into $M(\hat A)^\infty$.
Since $\phi^n=\hat\varepsilon$ on harmonic elements and
$\phi^\infty j_n=\phi^n$, we have
$\phi^\infty\theta=\hat\varepsilon$.

\smallskip

Izumi observed that $H^\infty(M(\hat A),P_\phi)$ can be given a
nicer description by embedding it into a larger algebra. Namely,
let $U$ be a unitary corepresentation such that the set of all its
irreducible components (irrespective of multiplicities) coincides
with $\supp\phi$. Then there exists a normal state~$\tilde\phi$
on~$B(H_U)$ such that $\tilde\phi|_{B_0(H_U)}$ is
$\alpha_U$-invariant and $\tilde\phi\pi_U=\phi$. Set
$(N,\tilde\phi^\infty)=\otimes_{-\infty}^{-1}(B(H_U),\tilde\phi)$.
The coactions $\alpha_{U^{\times n}}$ of $(A,\D)$ on
$\otimes^{-1}_{-n}B_0(H_U)$ define an ITP coaction $\alpha$ of
$(M,\D)$ on the von Neumann algebra $N$. Here $M=\pi_r(A)''$ is
the weak operator closure of $A$ in the regular
representation~$\pi_r$, so we are considering quantum groups and
their coactions in the von Neumann setting. The homomorphism
$\otimes^{-1}_{-\infty}\pi_U$ defines an embedding of $(M(\hat
A)^\infty,\phi^\infty)$ into $(N,\tilde\phi^\infty)$, and if we
identify $B(H_U\otimes\dots\otimes H_U)$ with
$\otimes^{-1}_{-n}B(H_U)\subset N$, then $j_n=\pi_{U^{\times n}}$.
Thus $j_n(M(\hat A))$ coincides with the relative commutant
$$
(B(H_U\otimes\dots\otimes H_U)^{\alpha_{U^{\times n}}})'\cap
B(H_U\otimes\dots\otimes H_U).
$$
Hence $\theta(H^\infty(M(\hat A),P_\phi))\subset N$
coincides with the relative commutant $(N^\alpha)'\cap N$.

\bigskip\bigskip

\section{Martin Boundary} \label{s3}

Throughout this section $\phi$ is presumed to be a generating
positive linear functional in $\C$ such that $\|\phi\|\leq 1$. If
$\phi$ is a state, we shall furthermore assume that it is
transient. Note that if $\|\phi\|\ <1$, then $M(\hat{A})\subset
D(G_\phi)$ and $\|{G_\phi}|_{M(\hat A)}\|=(1-\|\phi\|)^{-1}$.

\begin{defin}\label{3.1}
The Martin kernel for $P_\phi$ is the linear map $K_{\check\phi}
:\hat\A\to M(\hat A)$ given by
$$K_{\check\phi}(x)=G_{\check\phi}(x)G_{\check\phi}(I_0)^{-1}.$$
\end{defin}

According to Lemma \ref{2.4}(v) and the proof of Corollary
\ref{2.5}, the linear functional $\check\phi$ is generating and
transient. Thus $G_{\check\phi}(x)$ is a well-defined element of
$M(\hat A)$ for any $x\in\hat\A$. Since $\check\phi$ is
generating, the element $G_{\check\phi}(I_0)I_s$ is a non-zero
scalar multiple of $I_s$, so $G_{\check\phi}(I_0)$ is invertible
in $M(\hat\A)$ for any $s\in I$. A priori, therefore, we have
$K_{\check\phi}(x)\in M(\hat\A)$. However, observe that since
$\check\phi$ is generating, any positive element $x\in\hat\A$ is
majorized by a scalar multiple of
$\sum^n_{k=0}P^n_{\check\phi}(I_0)$ for some $n\in\7N$. As
$G_{\check\phi}P^n_{\check\phi}(I_0)\leq G_{\check\phi}(I_0)$, we
see that for any $x\in\hat\A_+$, there exists $c>0$ (depending on
$x$) such that $G_{\check\phi}(x)\leq cG_{\check\phi}(I_0)$. Hence
$K_{\check\phi} (\hat\A)\subset M(\hat A)$. As
$G_{\check\phi}(I_0)$ is central, obviously $K_{\check\phi}$ is
completely positive.

\begin{defin}\label{3.2}
The Martin compactification of the discrete quantum group $(\hat
A,\hat\D)$ with respect to~$P_\phi$ is the C$^*$-subalgebra
$\tilde A_\phi$ of $M(\hat A)$ generated by
$K_{\check\phi}(\hat\A)$ and $\hat A$. The Martin boundary
$A_\phi$ is the quotient C$^*$-algebra $\tilde{A}_\phi /\hat A$.
\end{defin}

We may think of $\hat A$ as the algebra of functions on our
discrete quantum group tending to zero at infinity, and of $M(\hat
A)$ as the algebra of all bounded functions. With this picture
in mind any
unital C$^*$-subalgebra of $M(\hat A)$ containing $\hat A$
plays the role of the algebra of continuous functions on some
compactification of the discrete quantum group. Note that $\tilde
A_\phi$ is unital since $K_{\check\phi}(I_0)=1$.

\smallskip

Suppose $(A,\D)=(C^*_r(\Gamma),\D)$, where $\Gamma$ is a discrete
group. A state $\phi$ on $\hat A=c_0(\Gamma)$ is represented by a
measure $\mu$. Then $p_\phi(h,g)=\mu(hg^{-1})$. Set
$G(h,g)=\sum_np_{\phi^n}(h,g)=\sum_n(\mu^{*n})(hg^{-1})$ and
$\displaystyle K(h,g)=\frac{G(h,g)}{G(e,g)}$, where $e=0$ is the
unit in $\Gamma$. We have $K_{\check\phi}(I_h)=\sum_gK(h,g)I_g$.
Thus $\tilde A_\phi=C(\bar\Gamma)$, where $\bar\Gamma$ is the
minimal compactification of $\Gamma$ for which all the functions
$g\mapsto K(h,g)$, $h\in\Gamma$, are continuous.

\subsection{Integral Representation of Superharmonic Elements}

In Subsection \ref{1.4} we introduced a sesquilinear form
$(\cdot,\cdot)_\eta$ associated to a KMS-state $\eta$. Here we use
a similar definition for the Haar weight $\hat\psi$, so
$(x,y)_{\hat\psi}=\hat\psi (x\sigma^{\hat\psi}_{-\ii}(y^*))$. We
need not worry about domain problems, since we shall always assume
that either $x$, or $y$ belongs to $\hat\A$.

\begin{thm} \label{3.3}
\mbox{\ }
\enu {i}
For any superharmonic element $x\in M(\hat\A)$, there exists a
bounded positive linear functional $\omega$ on $\tilde{A}_\phi$
such that $(y,x)_{\hat\psi}=\omega K_{\check\phi}(y)$ for
$y\in\hat\A$.
\enu {ii}
Conversely, for any bounded positive linear functional $\omega$ on
$\tilde{A}_\phi$ there exists a unique superharmonic element
$x_\omega$ such that $(y,x_\omega )_{\hat\psi} =\omega
K_{\check\phi}(y)$ for all $y\in\hat\A$. If $x_\omega$ is
harmonic, then $\omega|_{\hat A}=0$. Moreover, if $\supp\phi$ is
finite, then $x_\omega$ is harmonic if and only if $\omega|_{\hat
A}=0$.
\end{thm}

To prove the theorem, we need the following result.

\begin{lem} \label{3.4}
For $x,y\in\hat\A$ we have
$(P_\phi(x),y)_{\hat\psi}=(x,P_{\check\phi}(y))_{\hat\psi}$.
\end{lem}

\bp
It is enough to consider $\phi\in\C\cap\A$. Strong left
invariance of the Haar weight $\hat\varphi$ reads as
\begin{equation} \label{e3.1}
\hat\varphi(yP_\phi(x))=\hat\varphi(P_{\phi\hat S}(y)x).
\end{equation}
Recall that $\hat\psi=\hat\varphi(\rho^{-2}\,\cdot)$ on $\hat\A$.
Since $\hat\D(\rho^{-2})=\rho^{-2}\otimes\rho^{-2}$, we have
$P_{\phi\hat S}(\rho^{-2}y)=\rho^{-2}P_{\phi\hat
S(\rho^{-2}\,\cdot)}(y)$. We claim that $\phi\hat
S(\rho^{-2}\,\cdot)=\check\phi$. To prove this we may suppose
$\phi=\phi_s$. Since $\phi_s$ is $\sigma^{\hat\psi}$-invariant,
$\hat R$ is a $*$-antiisomorphism of $B(H_s)$ onto $B(H_{\bar
s})$, and $\hat R(\rho)=\rho^{-1}$, we have
$$
\phi_s\hat S=\phi_s\hat\tau_{-\ii}\hat R
=\phi_s\sigma^{\hat\psi}_\ii\hat R
=\phi_s\hat R=\frac{1}{d_s}\Tr\,\pi_{\bar s}(\rho\,\cdot),
$$
so $\phi_s\hat S(\rho^{-2}\,\cdot)=\phi_{\bar s}=\check\phi_s$.
Thus replacing $y$ by $\rho^{-2}y$ in (\ref{e3.1}) we get
$\hat\psi(yP_\phi(x))=\hat\psi(P_{\check\phi}(y)x)$. It remains to
replace $y$ by $\sigma^{\hat\psi}_\ii(y^*)$ and note that
$P_{\check\phi}\sigma^{\hat\psi}_t=\sigma^{\hat\psi}_tP_{\check\phi}$.
The latter equality follows from $\hat\D\sigma^{\hat\psi}_t
=(\sigma^{\hat\psi}_t\otimes\sigma^{\hat\psi}_t)\hat\D$ and the
fact that $\check\phi$ is $\sigma^{\hat\psi}$-invariant.
\ep

\bpp{Theorem \ref{3.3}}
To prove (i) first limit to the case when $x$ is a potential of an
element $x_0\in\hat\A_+$. Then we can take $\omega$ to be the
restriction of the positive linear functional
$(\cdot\,G_{\check\phi}(I_0),x_0)_{\hat\psi}$ to $\tilde A_\phi$.
Indeed, by the previous lemma we have
$$
\omega K_{\check\phi}(y)=(G_{\check\phi}(y),x_0)_{\hat\psi}
=(y,G_\phi(x_0))_{\hat\psi}=(y,x)_{\hat\psi}
$$
for any $y\in\hat\A$.
A general superharmonic element can be approximated from below by
potentials due to the balayage theorem. Hence it can be
approximated by potentials of elements in $\hat\A$. Thus there
exists a net $\{x_i\}_i$ of positive elements in $\hat\A$ such
that $G_\phi(x_i)\le x$ and $G_\phi(x_i)\to x$ (in the topology of
$M(\hat\A)$). Let $\omega_i$ be any positive linear functional on
$\tilde A_\phi$ satisfying
$\omega_iK_{\check\phi}(y)=(y,G_\phi(x_i))_{\hat\psi}$ for all
$y\in\hat\A$. Note that
$$
\omega_i(1)=\omega_iK_{\check\phi}(I_0)
=(I_0,G_\phi(x_i))_{\hat\psi} =\hat\varepsilon
G_\phi(x_i)\le\hat\varepsilon(x).
$$
Thus we can take for $\omega$ any weak$^*$ limit point of the net
$\{\omega_i\}_i$.

\smallskip

To prove (ii) note that the pairing $(\cdot,\cdot)_{\hat\psi}$ on
$\hat\A\times M(\hat\A)$ defines an antilinear order isomorphism
between $(\hat\A)'$ and $M(\hat\A)$. Hence, since $\omega
K_{\check\phi}$ is a positive linear functional on $\hat\A$, there
exists a unique positive element $x_{\omega}\in M(\hat\A)$ such
that $(\cdot,x_\omega)_{\hat\psi}=\omega K_{\check\phi}$. For any
$y\in\hat\A_+$ and $X\in\I$, we have
$$
G_{\check\phi}(y)=G_{\check\phi}P_{\check\phi}(y)+y
\ge G_{\check\phi}F_XP_{\check\phi}(y)+y,
$$
so $K_{\check\phi}(y)\ge K_{\check\phi}F_XP_{\check\phi}(y)
+yG_{\check\phi}(I_0)^{-1}$, whence
\begin{eqnarray*}
(y,P_\phi F_X(x_\omega))_{\hat\psi}
&=&(F_XP_{\check\phi}(y),x_\omega)_{\hat\psi}
=\omega K_{\check\phi}F_XP_{\check\phi}(y)
\le\omega K_{\check\phi}(y)-\omega(yG_{\check\phi}(I_0)^{-1})\\
&=&(y,x_\omega)_{\hat\psi}-\omega(yG_{\check\phi}(I_0)^{-1})
\le(y,x_\omega)_{\hat\psi}.
\end{eqnarray*}
Hence $P_\phi F_X(x_\omega)\le x_\omega$. It follows that
$x_\omega\in D(P_\phi)$ and $P_\phi(x_\omega)\le x_\omega$. If
$x_\omega$ is harmonic, then $(y,P_\phi
F_X(x_\omega))_{\hat\psi}\nearrow(y,P_\phi(x_\omega))_{\hat\psi}
=(y,x_\omega)_{\hat\psi}$, so we see that
$\omega(yG_{\check\phi}(I_0)^{-1})=0$, that is $\omega|_{\hat
A}=0$. Conversely, if $\omega|_{\hat A}=0$ and $\supp\phi$ is
finite, then $P_{\check\phi}(y)\in\hat\A$ and
$K_{\check\phi}(y)=K_{\check\phi}P_{\check\phi}(y)
+yG_{\check\phi}(I_0)^{-1}$, so
$(y,P_\phi(x_\omega))_{\hat\psi}=(y,x_\omega)_{\hat\psi}
-\omega(yG_{\check\phi}(I_0)^{-1})=(y,x_\omega)_{\hat\psi}$. Thus
$P_\phi(x_\omega)=x_\omega$.
\ep

\subsection{Canonical Coactions} \label{s3.2}

In Subsection \ref{s1.4} we introduced several coactions. The aim
of this subsection is to show that they induce coactions on the
Martin compactification and the Martin boundary.

\begin{thm} \label{3.5}
\mbox{\ }
\enu{i}
The right coaction $\hat\D$ of $(\hat A,\hat\D)$ on $\hat A$ has
the property $\hat\D(\tilde A_\phi)\subset M(\tilde
A_\phi\otimes\hat A)$, so it induces right coactions of $(\hat
A,\hat\D)$ on $\tilde A_\phi$ and $A_\phi$. We denote these
coactions by the same letter $\hat\D$.
\enu{ii}
The left coaction $\Phi$ of $(A,\D)$ on $\hat A$ has the property
$\Phi(\tilde A_\phi)\subset A\otimes\tilde A_\phi$. It induces
left coactions of $(A,\D)$ on $\tilde A_\phi$ and $A_\phi$, which
we again denote by $\Phi$.
\end{thm}

The proof is based on the following result.

\begin{lem}\cite[Lemma 2.2]{I} \label{3.6}
For any $x\in M(\hat A)$ we have
\enu{i}
$\hat\D P_\phi(x)=(P_\phi\otimes\iota)\hat\D(x)$;
\enu{ii}
$\Phi P_\phi(x)=(\iota\otimes P_\phi)\Phi(x)$.
\epp
\end{lem}

This lemma immediately implies part (ii) of Theorem \ref{3.5}.
Indeed, we have $\Phi G_{\check\phi}(x)=(\iota\otimes
G_{\check\phi})\Phi$ on $\hat\A$ (this expression makes sense,
since $\Phi(\hat\A)\subset\A\odot\hat\A$). As
$G_{\check\phi}(I_0)\in Z(M(\hat A))$, we have $\Phi
G_{\check\phi}(I_0)=1\otimes G_{\check\phi}(I_0)$, so
\begin{equation} \label{e3.2}
\Phi K_{\check\phi}(x)=(\iota\otimes K_{\check\phi})\Phi(x).
\end{equation}
This implies $\Phi(\tilde A_\phi)\subset A\otimes\tilde A_\phi$.
Furthermore, since $(\varepsilon\otimes\iota)\Phi=\iota$ on
$\hat\A$ (see Proposition \ref{1.7}), the assumptions of Corollary
\ref{1.4}(i) are fulfilled for the $*$-algebra generated by
$K_{\check\phi}(\hat\A)$ and $\hat\A$, so $\Phi\colon\tilde
A_\phi\to A\otimes \tilde A_\phi$ is a coaction.

\smallskip

Similarly, part (i) of Lemma \ref{3.6} implies $\hat\D
G_{\check\phi}(x)=(G_{\check\phi}\otimes\iota)\hat\D(x)$ for
$x\in\hat\A$. Here we consider $G_{\check\phi}\otimes\iota$ as an
operator going from the algebraic direct product $\prod_{s\in
I}\hat\A\otimes B(H_s)$ into $\prod_{s\in I}M(\hat\A)\otimes
B(H_s) =M(\hat\A\odot\hat\A)$. Thus
\begin{equation} \label{e3.3}
\hat\D K_{\check\phi}(x)
=(K_{\check\phi}\otimes\iota)\hat\D(x)(G_{\check\phi}(I_0)\otimes1)
(\hat\D G_{\check\phi}(I_0)^{-1})
=(K_{\check\phi}\otimes\iota)\hat\D(x)
((K_{\check\phi}\otimes\iota)\hat\D(I_0))^{-1}.
\end{equation}
Clearly $((K_{\check\phi}\otimes\iota)\hat\D(I_0))^{-1}\in
M(\hat\A\odot\hat\A)$. To show that $\hat\D K_{\check\phi}(x)$
belongs to the multiplier algebra of $\tilde A_\phi\otimes\hat A$,
first notice that $M(\tilde A_\phi\otimes\hat A)$ is the
C$^*$-algebraic direct product of $\tilde A_\phi\otimes B(H_s)$,
$s\in I$. Since $\hat\D K_{\check\phi}(x)\in M(\hat A\otimes\hat
A)$, we just have to show that
$(K_{\check\phi}\otimes\iota)((1\otimes I_s)\D(I_0))$ is
invertible in $M(\hat A)\otimes B(H_s)$, which entails $(1\otimes
I_s)((K_{\check\phi}\otimes\iota)\D(I_0))^{-1}\in\tilde
A_\phi\otimes B(H_s)$. To this end we need the following two
lemmas.

\begin{lem} \label{3.7}
Let $r,t\in I$ and suppose $p_\phi(r,t)\ne0$. Then
$I_rP_\phi\colon B(H_t)\to B(H_r)$ is a faithful completely
positive map.
\end{lem}

\bp
If $p_\phi(r,t)\ne0$, then $\phi_rP_\phi=\phi\phi_r$ majorizes a
scalar multiple of the state $\phi_t$, which obviously is faithful
on $B(H_t)$.
\ep

As a corollary, any non-zero superharmonic element is
invertible in $M(\hat\A)$.

\begin{lem} \label{3.8}
Consider $t\in I$ and a finite dimensional C$^*$-algebra $B$.
Suppose $x\in B(H_t)\otimes B$ is a positive element such that
$(\omega\otimes\iota)(x)$ is invertible for some state $\omega$ on
$B(H_t)$ (equivalently, the support of $x$ cannot be majorized by
a non-trivial projection $I_t\otimes p$). Then
$(K_{\check\phi}\otimes\iota)(x)$ is invertible in $\tilde
A_\phi\otimes B$.
\end{lem}

\bp
Pick $n\in\7N$ such that $p_{\check\phi^n}(0,t)\ne0$. By the
previous lemma $I_0 P_{\check\phi}^n$ is a faithful
completely positive map on $B(H_t)$, so
$P_{\check\phi}^n\ge c\omega(\cdot)I_0$ for some $c>0$, whence
$P_{\check\phi}^n\otimes\iota\ge c\omega(\cdot)I_0\otimes\iota$
on $B(H_t)\otimes B$. Then
$$
(G_{\check\phi}\otimes\iota)(x)
\ge (G_{\check\phi}P_{\check\phi}^n\otimes\iota)(x)
\ge c(G_{\check\phi}\otimes\iota)
    (\omega(\cdot)I_0\otimes\iota)(x)
=cG_{\check\phi}(I_0)\otimes((\omega\otimes\iota)(x)),
$$
so $(K_{\check\phi}\otimes\iota)(x)\ge
c1\otimes((\omega\otimes\iota)(x))$ is invertible.
\ep

In particular, the element $K_{\check\phi}(x)$ is invertible in
$\tilde A_\phi$ for any non-zero positive $x\in\hat\A$.

Recall that by Lemma \ref{1.2}, $(1\otimes I_s)\hat\D(I_0)\in
B(H_{\bar s})\otimes B(H_s)$ and $(\phi_{\bar
s}\otimes\iota)((1\otimes I_s)\hat\D(I_0))=\frac{1}{d^2_s}I_s$.
Therefore $(K_{\check\phi}\otimes\iota)((1\otimes
I_s)\hat\D(I_0))$ is invertible in $\tilde A_\phi\otimes B(H_s)$,
and thus $\hat\D(\tilde A_\phi)\subset M(\tilde A_\phi\otimes\hat
A)$. Consequently, $\hat\D\colon\tilde A_\phi\to M(\tilde
A_\phi\otimes\hat A)$ is a coaction by Corollary \ref{1.4}(ii).

\subsection{Poisson Boundary and States Representing the Unit}
\label{s3.3}

Throughout this subsection we assume in addition that $\phi$ is a
state, so the unit is a harmonic element. Let $\nu$ be a state on
$\tilde A_\phi$ representing the unit, that is $\nu
K_{\check\phi}=(\cdot,1)_{\hat\psi}=\hat\psi$ on $\hat\A$. Since
the unit is harmonic, the restriction of $\nu$ to $\hat A$ is zero
by Theorem \ref{3.3}(ii). Thus $\nu$ can be considered as a state
on the boundary $A_\phi$. Suppose in addition that $\nu$ is KMS.
Then, as in Subsection \ref{1.4}, we can define a dual map
$K^*_{\check\phi}\colon\pi_\nu(A_\phi)''\to M(\hat\A)$
determined by
$(K_{\check\phi}(x),a)_\nu=(x,K^*_{\check\phi}(a))_{\hat\psi}$ for
all $x\in\hat\A$ and $a\in\pi_\nu(A_\phi)''$. As
$(x,K^*_{\check\phi}(1))_{\hat\psi}=\nu
K_{\check\phi}(x)=\hat\psi(x)$, we deduce that
$K^*_{\check\phi}(1)=1$. Thus~$K^*_{\check\phi}$ is, in fact, a
normal unital completely positive map from $\pi_\nu(A_\phi)''$
into $M(\hat A)$. Equivalently, it can be described as follows: if
$a\in\pi_\nu(A_\phi)''$ is positive, then $K^*_{\check\phi}(a)$ is
a superharmonic element corresponding to the positive linear
functional $(\cdot,a)_\nu$.

\begin{lem} \label{3.9}
We have $\Im\, K^*_{\check\phi}\subset H^\infty(M(\hat
A),P_\phi)$.
\end{lem}

\bp
Since $(\cdot,a)_\nu\le\nu(a)\nu$ for $a\ge0$, it suffices to
prove the following: if $\omega_1$ and $\omega_2$ are positive
linear functionals on $\tilde A_\phi$ such that
$\omega_1\le\omega_2$ and $x_{\omega_2}$ is harmonic, then
$x_{\omega_1}$ is also harmonic. Since $x_{\omega_2}$ is
harmonic and both $x_{\omega_1}$ and
$x_{\omega_2}-x_{\omega_1}=x_{\omega_2-\omega_1}$ are
superharmonic, the element $x_{\omega_1}$ must be harmonic.
\ep

Notice also that as $K_{\check\phi}(I_0)=1$, we have
$\hat\varepsilon K_{\check\phi}^*=\nu$.

\smallskip

Define a dynamics $\gamma$ on $\tilde A_\phi$ as the restriction
of $\sigma^{\hat\psi}$ to $\tilde A_\phi$. As already remarked in
the proof of Lemma \ref{3.4}, we have the equality
$\sigma^{\hat\psi}_t
P_{\check\phi}=P_{\check\phi}\sigma^{\hat\psi}_t$, so
$\sigma^{\hat\psi}_t
K_{\check\phi}=K_{\check\phi}\sigma^{\hat\psi}_t$. Hence $\gamma$
is a strongly continuous one-parameter automorphism group on
$\tilde A_\phi$. Since $\gamma_t(\hat A)=\hat A$, we also have a
dynamics on $A_\phi$, which we denote by the same letter $\gamma$.
For the coactions $\Phi$ and $\hat\D$ on $\tilde A_\phi$ and
$A_\phi$ we have
\enu{i}
$\Phi\gamma_t=(\tau_t\otimes\gamma_t)\Phi$;
\enu{ii}
$\hat\D\gamma_t=(\gamma_t\otimes\hat\tau_{-t})\hat\D$.

\begin{thm} \label{3.10}
Let $\nu$ be a weak$^*$ limit point of the sequence
$\{\phi^n|_{\tilde A_\phi}\}^\infty_{n=1}$. Then $\nu$ is
a $\gamma$-KMS state representing the unit.
\end{thm}

\bp
Since $\phi^n|_{\tilde A_\phi}$ is a $\gamma$-KMS state,
clearly $\nu$ is $\gamma$-KMS.

Consider the positive linear functional $\nu K_{\check\phi}$ on
$\hat\A$. We assert that it is $\sigma^{\hat\psi}$-KMS. Since
$\nu$ is a weak$^*$ limit point of $\phi^n$, and $\phi^n$ is a
convex combination of the states $\phi_s$, it is enough to check
that $\phi_s K_{\check\phi}$ is $\sigma^{\hat\psi}$-KMS for any
$s\in I$. As $G_{\check\phi}(I_0)$ is central, the element
$I_sG_{\check\phi}(I_0)^{-1}$ is a scalar multiple of~$I_s$, so we
need only to consider the linear functional
$\phi_sG_{\check\phi}$. But it is the sum of
$\sigma^{\hat\psi}$-KMS functionals $\check\phi^n\phi_s$, so it
must be $\sigma^{\hat\psi}$-KMS itself. Thus our assertion is
proved.

Since $\phi_s$ is a unique (up to a scalar)
$\sigma^{\hat\psi}$-KMS functional on $B(H_s)$, to prove that
$\nu$ represents the unit, that is $\nu K_{\check\phi}=\hat\psi$ on
$\hat\A$, it is enough to verify this equality on $Z(\hat\A)$. In
this case we can apply the classical theory \cite[Theorem
7.2.7]{R}. Indeed, Lemma \ref{3.4} (or Lemma \ref{2.4}(v)) shows
that the classical Markov operators $P_\phi|_{Z(\hat A)}$ and
$P_{\check\phi}|_{Z(\hat A)}$ are in duality with respect to
$\hat\psi|_{Z(\hat A)}$. Therefore $K_{\check\phi}|_{Z(\hat\A)}$
coincides with the operator $\hat K$ introduced in
\cite[Proposition 7.2.3]{R}. Then with the notation of our
Subsection \ref{s2.1}, Theorem 7.2.7 in \cite{R} states that for
any $x\in Z(\hat\A)$ the sequence
$\{j_nK_{\check\phi}(x)\}^\infty_{n=1}\subset
L^\infty(\Omega,\7P_0)\subset M(\hat A)^\infty$ converges a.e. to
an element $j_\infty K_{\check\phi}(x)$ (in particular, it
converges in strong$^*$ operator topology on $M(\hat A)^\infty$)
and $\phi^\infty j_\infty K_{\check\phi}=\hat\psi$ on $Z(\hat\A)$.
Since $\phi^\infty j_n=\phi^n$, this is the same as saying that $\nu
K_{\check\phi}=\hat\psi$ on $Z(\hat\A)$.
\ep

It is shown in \cite{I} that if we consider our quantum groups in
the von Neumann setting, then $\hat\D$ and $\Phi$ define right and
left coactions of $(\hat M,\hat\D)$ and $(M,\D)$ on
$H^\infty(M(\hat A),P_\phi)$, respectively. Here $\hat M=M(\hat
A)$ and $M$ denote the weak operator closures of $\hat A$ and $A$,
respectively, in their regular representations.

\begin{prop} \label{3.11}
Let $\nu$ be a $\gamma$-KMS state on $A_\phi$ representing the
unit. Regard $K^*_{\check\phi}$ as a map from $A_\phi$ to
$H^\infty(M(\hat A),P_\phi)$. Then
\enu{i}
$K^*_{\check\phi}$ intertwines the left coactions of $(A,\D)$ if
and only if $\nu$ is invariant;
\enu{ii}
$K^*_{\check\phi}$ intertwines the right coactions of $(\hat
A,\hat\D)$ if and only if $\nu$ is quasi-invariant with
Radon-Nikodym cocycle $(K_{\check\phi}\otimes\iota)\hat\D(I_0)$;
moreover, if these two equivalent conditions are satisfied, then
$K^*_{\check\phi}=(\nu\otimes\iota)\hat\D$.
\end{prop}

\bp
To prove (i) consider for any $\omega\in A^*$ the operator
$S_\omega$ on $M(\hat A)$ given by
$S_\omega=(\omega\otimes\iota)\Phi$. Because
$\hat\psi|_{B(H_s)}=d_s^2\phi_s$ is $\Phi$-invariant for any $s\in
I$, by Proposition \ref{1.6} we see that $S^*_\omega=S_{\bar\omega
R}$ with respect to the inner product $(\cdot,\cdot)_{\hat\psi}$.
Thus $(S_\omega(x),y)_{\hat\psi}=(x,S_{\bar\omega
R}(y))_{\hat\psi}$ if $x$ or $y$ is in $\hat\A$. Analogously,
consider the operator $T_\omega=(\omega\otimes\iota)\Phi$ on
$A_\phi$. Then $K^*_{\check\phi}$ intertwines the left coactions,
that is $\Phi K^*_{\check\phi}=(\iota\otimes
K^*_{\check\phi})\Phi$ on $A_\phi$, if and only if $S_\omega
K^*_{\check\phi}=K^*_{\check\phi} T_\omega$  for all $\omega\in
A^*$. On the other hand, the state $\nu$ is invariant if and only
if $\nu T_\omega=\omega(1)\nu$ for all $\omega\in A^*$. The
identity (\ref{e3.2}) can be rewritten as $T_\omega
K_{\check\phi}=K_{\check\phi}S_\omega$. For $a\in A_\phi$ and
$x\in\hat\A$, we have
$$
(S_\omega K^*_{\check\phi}(a),x)_{\hat\psi}
=(K^*_{\check\phi}(a), S_{\bar\omega R}(x))_{\hat\psi}
=(a,K_{\check\phi}S_{\bar\omega R}(x))_\nu
=(a,T_{\bar\omega R}K_{\check\phi}(x))_\nu
$$
and $(K^*_{\check\phi}T_\omega(a),x)_{\hat\psi}
=(T_\omega(a),K_{\check\phi}(x))_\nu$. So if $S_\omega
K^*_{\check\phi}=K^*_{\check\phi} T_\omega$, then taking $x=I_0$,
using $K_{\check\phi}(I_0)=1$ and
$T_{\bar\omega R}(1)=\overline{\omega(1)}1$, we obtain $\nu
T_\omega(a)=\omega(1)\nu(a)$. Conversely, if $\nu$ is invariant,
then $T^*_\omega=T_{\bar\omega R}$ by Proposition
\ref{1.6}, and we get $S_\omega K^*_{\check\phi}=
K^*_{\check\phi}T_\omega$.

\smallskip

To settle (ii) first notice that
$y=(K_{\check\phi}\otimes\iota)\hat\D(I_0)$ is a cocycle, i.e.
$(\iota\otimes\hat\D)(y)=(\hat\D\otimes\iota)(y)(y\otimes1)$. This
follows from coassociativity of $\hat\D$ and (\ref{e3.3}), which
can be written as $(K_{\check\phi}\otimes\iota)\hat\D(x)=\hat\D
K_{\check\phi}(x)y$. To see this, we compute
\begin{eqnarray*}
(\iota\otimes\hat\D)(y)
&=&(K_{\check\phi}\otimes\iota\otimes\iota)
    (\iota\otimes\hat\D)\hat\D(I_0)
    =((K_{\check\phi}\otimes\iota)\hat\D\otimes\iota)\hat\D(I_0)\\
&=&(\hat\D K_{\check\phi}\otimes\iota)\hat\D(I_0)(y\otimes1)
    =(\hat\D\otimes\iota)(y)(y\otimes1).
\end{eqnarray*}
Now for $\omega\in\A$, consider the operators
$S_\omega=(\iota\otimes\omega)\hat\D$ and
$T_\omega=(\iota\otimes\omega)\hat\D$ on $M(\hat A)$ and $A_\phi$,
respectively. As in Proposition \ref{1.6} strong right invariance
of the Haar weight $\hat\psi$ can be expressed by the equality
$S^*_\omega=S_{\bar\omega\hat R}$.

Suppose $K^*_{\check\phi}$ intertwines the coactions of $(\hat
A,\hat\D)$. Then for any $a\in A_\phi$ and $\omega\in\A$, we get
\begin{eqnarray*}
\nu T_\omega(a)
&=&\hat\varepsilon K^*_{\check\phi}T_\omega(a)
    =(K^*_{\check\phi}T_\omega(a),I_0)_{\hat\psi}
    =(S_\omega K^*_{\check\phi}(a),I_0)_{\hat\psi}\\
&=&(K^*_{\check\phi}(a),S_{\bar\omega\hat R}(I_0))_{\hat\psi}
    =(a,K_{\check\phi}S_{\bar\omega\hat R}(I_0))_\nu.
\end{eqnarray*}
Since $\gamma_t K_{\check\phi}=K_{\check\phi}\sigma^{\hat\psi}_t$
and $\hat\D\sigma^{\hat\psi}_t
=(\sigma^{\hat\psi}_t\otimes\hat\tau_{-t})\hat\D$, this equality
is equivalent to property (i) in Definition~\ref{1.9}. Thus $\nu$
is quasi-invariant with Radon-Nikodym cocycle
$y=(K_{\check\phi}\otimes\iota)\hat\D(I_0)$. Moreover, as
$(S_\omega K^*_{\check\phi}(a),I_0)_{\hat\psi} =\hat\varepsilon
S_\omega K^*_{\check\phi}(a) =\omega K^*_{\check\phi}(a)$ and $\nu
T_\omega(a)=\omega(\nu\otimes\iota)\hat\D(a)$, we have
$K^*_{\check\phi}(a)=(\nu\otimes\iota)\hat\D(a)$. Note that the
last argument does not use the definition of $K^*_{\check\phi}$.
The very fact that $K^*_{\check\phi}$ intertwines the coactions
implies $K^*_{\check\phi}
=(\hat\varepsilon K^*_{\check\phi}\otimes\iota)\hat\D$.

Conversely, suppose $\nu$ is quasi-invariant with Radon-Nikodym
cocycle $y=(K_{\check\phi}\otimes\iota)\hat\D(I_0)$. For any
$a,b\in A_\phi$ and $\omega\in\A$, Proposition
\ref{1.10}(ii) yields
$$
(\nu\otimes\omega)(\hat\D(a)(b\otimes1))
=(\nu\otimes\omega\hat S)((a\otimes 1)\hat\D(b)y).
$$
Taking $b=K_{\check\phi}(x)$ for $x\in\hat\A$, and using
$(K_{\check\phi}\otimes\iota)\hat\D(x)=\hat\D K_{\check\phi}(x)y$,
we obtain
$$
(\nu\otimes\omega)(\hat\D(a)(K_{\check\phi}(x)\otimes1))
=(\nu\otimes\omega\hat S)((a\otimes 1)
 (K_{\check\phi}\otimes\iota)\hat\D(x)).
$$
Replacing $x$ by $\sigma^{\hat\psi}_{-\ii}(x^*)$, this identity
can be rewritten as $(T_\omega(a),K_{\check\phi}(x))_\nu
=(a,K_{\check\phi}S_{\bar\omega\hat R}(x))_\nu$. It follows
that $K^*_{\check\phi}T_\omega=S_\omega K^*_{\check\phi}$ for all
$\omega\in\A$. In other words, $K^*_{\check\phi}$ intertwines the
right coactions.
\ep

If $\nu$ is an invariant state for the coaction $\Phi$ of $(A,\D)$
on $A_\phi$, it induces a coaction of $(A,\D)$ on~$\pi_\nu(A_\phi)$.
The latter can be extended to a coaction of
$(M,\D)$ on $\pi_\nu(A_\phi)''$. The reason for this is that there
exists a unitary on $H_{\varphi}\otimes H_\nu$ implementing this
coaction (see the proof of Proposition~\ref{1.5}).

If $\nu$ is a quasi-invariant $\gamma$-KMS state with
Radon-Nikodym cocycle $y=(K_{\check\phi}\otimes\iota)\hat\D(I_0)$
for the right coaction $\hat\D$ of $(\hat A,\hat\D)$ on $A_\phi$,
then this coaction is implemented on $H_\nu\otimes H_\varphi$ by
the unitary~$U$ given by
$$
U(J_\nu\otimes J_\varphi)y^\2(J_\nu\otimes J_\varphi)
  (a\xi_\nu\otimes\xi)
=\hat\D(a)(\xi_\nu\otimes\xi)\ \ \ \hbox{for}\ a\in A_\phi\ \
\hbox{and}\ \xi\in \A\xi_{\varphi}.
$$
To verify that $U$ is unitary, recall that $\hat R(x)=J_\varphi
x^* J_\varphi$ for $x\in M(\hat A)$. So if we denote by
$\omega_\xi$ the linear functional $(\cdot\,\xi,\xi)$ on $M(\hat
A)$, we get
\begin{eqnarray*}
\|(J_\nu\otimes J_\varphi)y^\2(J_\nu\otimes J_\varphi)
     (a\xi_\nu\otimes\xi)\|^2
&=&(\nu\otimes\omega_\xi\hat R)((a^*a\otimes1)
       (\gamma_{-\ii}\otimes\iota)(y))\\
&=&(\nu\otimes\omega_\xi\hat S)((a^*a\otimes 1)y)\\
&=&(\nu\otimes\omega_\xi)\hat\D(a^*a)
    =\|\hat\D(a)(\xi_\nu\otimes\xi)\|^2.
\end{eqnarray*}
It follows that $\hat\D$ induces a coaction of $(\hat A,\hat\D)$
on $\pi_\nu(A_\phi)$, which can be extended to a coaction of
$(\hat M,\hat\D)$ on $\pi_\nu(A_\phi)''$. (In fact, Proposition
\ref{1.11} shows that the particular form of the cocycle $y$ does
not play any role in this argument.)

\smallskip

A detailed study of connections between Martin and Poisson
boundaries will be given in a subsequent paper. For the
computation in the next section the following result is
sufficient.

\begin{prop} \label{3.12}
Suppose that the Martin kernel $K_{\check\phi}$ considered as a
map from $\hat\A$ into $A_\phi$ has dense range. Then
\enu{i}
the map $\omega\mapsto x_\omega$, which associates a superharmonic
element to a positive linear functional on~$A_\phi$, is injective;
\enu{ii}
if $\nu$ is the unique state on $A_\phi$ representing the unit,
then $\nu$ is $\gamma$-KMS and the map
$K_{\check\phi}^*\colon\pi_\nu(A_\phi)''\to H^\infty(M(\hat
A),P_\phi)$ is an isomorphism of von Neumann algebras
intertwining the left coactions of $(M,\D)$.
\end{prop}

\bp
Part (i) is obvious as $x_\omega$ by definition determines the
values of $\omega$ on $K_{\check\phi}(\hat\A)$. It remains to
prove (ii). Since there always exists a $\gamma$-KMS state
representing the unit by Theorem~\ref{3.10}, the state $\nu$ must
be $\gamma$-KMS and, being regarded as a state on $\tilde A_\phi$,
be the weak$^*$ limit of the sequence $\{\phi^n|_{\tilde
A_\phi}\}^\infty_{n=1}$. As the states $\phi^n$ are invariant with
respect to the left coaction of $(\hat A,\hat\D)$ on $\tilde
A_\phi$, the state $\nu$ is invariant as well. (Another way to see
this is to note that the state $(\varphi\otimes\nu)\Phi$ also
represents the unit and is invariant.) It follows that we have a
well-defined coaction of $(M,\D)$ on $\pi_\nu(A_\phi)''$, and
$K_{\check\phi}^*$ intertwines this coaction with the coaction on
$H^\infty(M(\hat A),P_\phi)$. We need now only to prove that
$K_{\check\phi}^*$ is an isomorphism. Since
$K_{\check\phi}(\hat\A)$ is dense in $\pi_\nu(A_\phi)''$,
obviously $K_{\check\phi}^*$ is injective. Let $x\in
H^\infty(M(\hat A),P_\phi)$ be positive and non-zero. Suppose
$\omega$ is the unique positive linear functional on $A_\phi$
representing $x$. Then the linear functional
$\nu-\|x\|^{-1}\omega$ represents the positive harmonic element
$1-\|x\|^{-1}x$. Since such a functional is unique and positive,
clearly $\omega\le\|x\|\nu$. Hence $\omega=(\cdot,a)_\nu$ for some
positive $a\in\pi_\nu(A_\phi)''$. Then $K_{\check\phi}^*(a)=x$, so
$K_{\check\phi}^*$ is a bijection. Moreover, it maps positive
elements of $\pi_\nu(A_\phi)''$ {\it onto} positive elements of
$H^\infty(M(\hat A),P_\phi)$, so the inverse map is positive.
Since $K_{\check\phi}^*$ is also unital and $2$-positive, it is a 
$*$-isomorphism (see e.g. \cite[Corollaries~2.2 and~3.2]{C}).
\ep

\bigskip\bigskip

\section{Martin Boundary of the Dual of Quantum SU(2)} \label{s4}

In this section $(A,\D)$ will denote the compact quantum group
$SU_q (2)$ of Woronowicz \cite{Wor}. We assume that the
deformation parameter $q$ lies in $(0,1)$, but our results are
also true for $q\in (-1,0)$. So $A$ is the universal unital
C$^*$-algebra with generators $\alpha$ and $\gamma$ satisfying the
relations
$$
\alpha^*\alpha+\gamma^*\gamma =1,
\ \ \alpha\alpha^* +q^2\gamma^*\gamma =1,\ \
\gamma^*\gamma =\gamma\gamma^* ,
$$
$$
\alpha\gamma =q\gamma\alpha,\ \ \alpha\gamma^*=q\gamma^*\alpha .
$$
The comultiplication $\D$ is determined by the formulas
$$
\D (\alpha )=\alpha\otimes\alpha -q\gamma^*\otimes\gamma,\ \
\D (\gamma )=\gamma\otimes\alpha +\alpha^*\otimes\gamma .
$$
Recall that the Haar state $\varphi$ of $(A,\D )$ is faithful
\cite{N} (see also \cite{BMT}). The characters $f_{it}$,
$t\in\7R$, on $\A$ can be extended to bounded characters on $A$
which we denote by the same symbols.

Consider the C$^*$-subalgebra $B$ of $A$ given by
$$
B=\{a\in A\ |\ f_{it}* a =a\ \forall t\in\7R\}.
$$
Then $\D$ and $\hat\Phi$ (see Proposition \ref{1.7}) induce a left
coaction of $(A,\D)$ on $B$ and a right coaction of $(\hat
A,\hat\D)$ on $B$, respectively. The C$^*$-algebra $B$ with the
left coaction of $(A,\D)$ is the quantum homogeneous sphere of
Podle\'{s}. Note that Podle\'{s} considers right coactions, so he
deals with $\7T\backslash SU_q (2)$, while we consider $SU_q
(2)/\7T$.

We are now in the position to state our main results in this
section.

\begin{thm}\label{4.1}
Let $\phi=\sum_s\lambda_s\phi_s\in\C$ be a generating state.
Suppose $\sum_s\lambda_s\dim H_s<\infty$. Then there
exists an isomorphism $A_\phi\iso B$ which intertwines the left
coactions of $(A,\D)$ and the right coactions of $(\hat
A,\hat\D)$.
\end{thm}

Thus the Martin boundary of the dual of quantum $SU(2)$ is
identified with the quantum homogeneous $2$-sphere of Podle\'{s}.
The result is also valid for $\|\phi\|<1$. Under this assumption
the corresponding result for ordinary $SU(2)$ was first
established by Biane \cite{B4}.

\begin{thm}\label{4.2}
Retaining the assumptions on $\phi$ from Theorem \ref{4.1}, there
exists a unique state $\nu$ on $A_\phi$ representing the unit. The
dual map $K^*_{\check\phi}\colon\pi_\nu (A_\phi)''\to H^\infty
(M(\hat A), P_\phi)$ is an isomorphism intertwining the left
coactions of $(M,\D)$ and the right coactions of $(\hat M,\hat\D)$.
The composition of $K^*_{\check\phi}$ with the isomorphism $B\iso
A_\phi$ from Theorem \ref{4.1} is given by $b\mapsto
(\varphi\otimes\iota)\hat\Phi (b)$.
\end{thm}

Thus one can think of $H^\infty (M(\hat A), P_\phi)$ as the algebra
of bounded measurable functions on the quantum homogeneous sphere.

In the case when $\supp\phi$ is finite, the fact that the map
$\Theta =(\varphi\otimes\iota)\hat\Phi$ gives an isomorphism
of~$\pi_\varphi (B)''$ onto $H^\infty (M(\hat A), P_\phi)$ was
established by Izumi \cite{I}.

\smallskip

The rest of this section is devoted to formulating and proving
more precise versions of these results.

\subsection{Quantum Spheres and the Dual of Quantum SU(2)}

First recall that $\A$ is the $*$-subalgebra of $A$ generated by
$\alpha$ and $\gamma$. The coinverse and the counit on $\A$ are
given by
$$
S(\alpha )=\alpha^*,\ \ S(\gamma )=-q\gamma,
\ \ \varepsilon (\alpha )=1, \ \ \varepsilon(\gamma)=0.
$$
The characters $f_z$ are defined by
\begin{equation}\label{e4.1}
f_z (\gamma)=0,\ \ f_z(\alpha)=q^{-z}.
\end{equation}

\smallskip

From the point of view of representation theory it is convenient
to introduce the quantized universal enveloping algebra $\sl2$. It
is by definition the universal unital $*$-algebra generated by
elements $e,f,k,k^{-1}$ satisfying the relations
$$
kk^{-1}=k^{-1}k=1,\ \ ke=qek,\ \ kf=q^{-1}fk,\ \
ef-fe=\frac{1}{q-q^{-1}}(k^2-k^{-2}),
$$
$$
k^*=k,\ \ e^* =f.
$$
The algebra $\sl2$ is, in fact, a dense subalgebra of
$\A'=M(\hat\A)$ and the restrictions of $\hat\D$, $\hat S$ and
$\hat\varepsilon$ to $\sl2$ turn it into a Hopf $*$-algebra.
Explicitly,
$$
\hat\D (k)=k\otimes k,\ \
\hat\D (e)=e\otimes k^{-1}+k\otimes e,\ \
\hat\D (f)=f\otimes k^{-1}+k\otimes f,
$$
$$
\hat S(k)=k^{-1},\ \ \hat S(e)=-q^{-1}e,\ \ \hat S(f)=-qf,
$$
$$
\hat\varepsilon (k)=1,\ \
\hat\varepsilon (e)=\hat\varepsilon (f)=0.
$$

The set $I$ of equivalence classes of irreducible representations
is identified with the set $\2\7Z_+$ of non-negative
half-integers. Note that $\dim H_s =2s+1$. The basis
$\{\xi^s_i\}_{i=-s}^s$ for $H_s$ is chosen in such a way that
\begin{eqnarray}
\pi_s (k)\xi^s_i&=&q^{-i}\xi^s_i,\label{e4.3}\\
\pi_s (e)\xi^s_i&=&([s+i]_q[s-i+1]_q)^\2\xi^s_{i-1},\label{e4.4}\\
\pi_s (f)\xi^s_i&=&([s-i]_q[s+i+1]_q)^\2\xi^s_{i+1},\label{e4.5}
\end{eqnarray}
where $[n]_q$ is the $q$-number for $n$, so
$$
[n]_q =\frac{q^n -q^{-n}}{q-q^{-1}}.
$$
The fundamental corepresentation corresponds to spin $s=\2$.
Thus
\begin{equation}\label{e4.2}
U^\2 =(u^\2_{ij})_{ij}=
\pmatrix{\alpha & -q\gamma^* \cr
\gamma & \alpha^*},
\end{equation}
so e.g. $u^\2_{-\2,-\2}=\alpha$ and $u^\2_{\2,-\2}=\gamma$. The
formulas (\ref{e4.3}--\ref{e4.5}) for $s=\2$ and (\ref{e4.2})
determine the pairing between $(\A,\D)$ and $(\sl2,\hat\D)$
uniquely. Since $k$ is a character on $\A$, the formulas
(\ref{e4.1}) and (\ref{e4.2}) show that
\begin{equation}\label{e4.6}
\rho=k^{-2}.
\end{equation}
In particular, for the quantum dimension we have $d_s =[2s+1]_q$.
Also the basis $\{\xi^s_i\}_i$ satisfies the conventions of
Subsection \ref{s1.3} in that $\pi_s (\rho)$ is diagonal.
Introducing matrix units $m^s_{ij}$ and $n^{\bar s}_{ij}$ as in
that subsection, we get the following identities.

\begin{lem}\label{4.3}
\mbox{\ }
\enu{i}
$s=\bar s$ and $n_{ij}^s =(-1)^{i-j}m^s_{-i,-j}$;
\vspace{2mm}
\enu{ii}
$\F (u^s_{ij})=(-1)^{i-j}d_s^{-1}q^{-i-j}m^s_{-i,-j}$ and
$(u^s_{ij})^*=(-1)^{i-j}q^{-i+j}u^s_{-i,-j}$;
\vspace{2mm}
\enu{iii}
$\alpha^* u^s_{jj}=d_s^{-1}(q^{-s+j}[s+j+1]_q
u^{s+\2}_{j+\2,j+\2}+q^{s+j+1}[s-j]_q u^{s-\2}_{j+\2,j+\2})$;
\vspace{2mm}
\enu{iv}
$\displaystyle u^s_{sj}=\left[2s\atop
s+j\right]^\2_{q^2}(\alpha^*)^{s+j}\gamma^{s-j}$.
\end{lem}

Here $\displaystyle\left[n\atop m\right]_r
=\frac{(r;r)_n}{(r;r)_m(r;r)_{n-m}}$, where $\displaystyle
(a;r)_n=\prod_{i=0}^{n-1}(1-ar^i)$ for $n\ge1$ and $(a;r)_0=1$.

\smallskip

\bpp{Lemma \ref{4.3}}
It is well-known that any unitary corepresentation of $SU_q (2)$
is self-conjugate. Explicitly, using (\ref{e4.3}--\ref{e4.5}) and
(\ref{e4.6}) it is easy to check that the unitary $\tilde
J_s\colon H_s\to\overline {H}_s$ given by $\tilde J_s\xi^s_j
=(-1)^{[j]}\overline{\xi^s_{-j}}$, where $[j]$ is the integral
part of $j$, intertwines $\pi_s$ and $\pi_{\overline{U^s}}$, as
$$
\pi_{\overline{U^s}}(x)\bar\xi
=\overline{\pi_s(\rho^{-\2}\hat S (x^*)\rho^\2)\xi}.
$$
Thus $s=\bar s$ and $\hat R (x)=J_s x^* J_s^{-1}$ for $x\in
B(H_s)$, where $J_s$ is the antilinear isometry on $H_s$ given by
$J_s\xi^s_j =(-1)^{[j]}\xi^s_{-j}$ (see Subsection \ref{1.3}).
Hence $n_{ij}^s =(-1)^{i-j}m^s_{-i,-j}$.
\smallskip

Now (ii) follows from (i) and Lemma \ref{1.1}. Part (iii) is a
particular case of the formulas for the Clebsch-Gordan
coefficients \cite{V} (see e.g. \cite{GKK} for more details). The
formula in (iv) is from \cite[Theorem~1.8]{MMNN}.
\ep

\medskip

Let $\alpha\colon B\to A\otimes B$ be a left coaction of $(A,\D)$
on a unital C$^*$-algebra $B$. For $s\in I$, consider the spectral
subspace
$$
B(s)=\span\{(\varphi\otimes\iota)((a^*\otimes1)\alpha(b))\ |\
b\in B,\ a\in\span\{u^s_{ij}\}_{i,j}\}.
$$
Podle\'{s} \cite{P} classified the coactions satisfying three
additional properties: $B(0)=\7C 1$, $B(1)$ has multiplicity one
(equivalently, $\dim B(1)=\dim H_1=3$), and $B$ is generated as a
C$^*$-algebra\linebreak by~$B(1)$. Such pairs $(B,\alpha)$ he
called quantum spheres. They are classified by a parameter
$c\in\{-(q^n+q^{-n})^{-2}\ |\ n=2,3,\dots\}\cup[0,\infty]$. The
corresponding algebra $B$ is denoted by $C(S^2_{q,c})$. We are not
interested in $C(S^2_{q,\infty})$. For $c\ne\infty$, the algebra
$C(S^2_{q,c})$ is the universal unital C$^*$-algebra generated by
elements $X_{-1}$, $X_0$ and $X_1$ satisfying the relations
$$
X^*_{-1}=-qX_1,\ \ X^*_0 =X_0,
$$
$$
X^2_0 +X^*_{-1}X_{-1}+X^*_1 X_1=1+(q+q^{-1})^2 c,
$$
\begin{equation}\label{e4.7}
qX_1X_0-q^{-1}X_0 X_1=(q^{-1}-q)X_1,
\end{equation}
\begin{equation}\label{e4.8}
(q^{-1}-q)X^2_0+X_{-1}X_1-X_1 X_{-1}=-(q^{-1}-q)X_0.
\end{equation}
The coaction $\alpha$ is uniquely defined by the identity
$$
\alpha (X_j)=\sum_{k=-1}^1 u^1_{jk}\otimes X_k ,\ \
j=-1,0,1.
$$
For $c=-(q^{2s+1}+q^{-2s-1})^{-2}$, $s\in\2\7N$, the quantum
sphere $(C(S^2_{q,c}),\alpha)$ is nothing else than
$(B(H_s),\alpha_s)$. Thus there exist uniquely determined
generators $X^s_j\in B(H_s)$, $j=-1,0,1$, to be displayed in the
result below.

\begin{lem}\label{4.4}
Set $\chi_{-1}=-qfk$,
$\displaystyle\chi_0=\frac{ef-q^2fe}{\sqrt{[2]_q}}$ and $\chi_1
=qek$. Then $X^s_j=\lambda^{-1}_s\pi_s (\chi_j)$, where
$\displaystyle\lambda_s
=\frac{q(q^{2s+1}+q^{-2s-1})}{(q-q^{-1})\sqrt{[2]_q}}$.
\end{lem}

\bp
See Lemma 4.1 and Lemma 5.7 in \cite{I}. For the sake of
convenience we sketch a proof.

The fact that $\tilde X^s_j=\pi_s (\chi_j)$, $j=-1,0,1$, is a
basis in the spectral subspace of $B(H_s)$ corresponding to spin
$1$, that is
\begin{equation}\label{e4.10}
\Phi(\tilde X^s_j)=\sum_{k=-1}^1 u^1_{jk}\otimes \tilde X^s_k ,
\end{equation}
is deduced using the adjoint action. Namely, if we set $\ad
X=(X\otimes\iota)\Phi$, the mapping $X\mapsto \ad X$ is an
antirepresentation (i.e. linear but antimultiplicative) of $\sl2$
on $\hat\A$. Thus, to check (\ref{e4.10}) it suffices to verify
the equality
$$
(\ad X)(\tilde X^s_j)
=\sum_{k=-1}^1 \pi_1(X)_{jk}\otimes \tilde X^s_k
$$
for $X=e,f,k,k^{-1}$. The coefficients $\pi_1(X)_{jk}$ can be read
off (\ref{e4.3}--\ref{e4.5}). On the other hand, the adjoint
action can be computed explicitly. Namely, if $\hat\D (X)=\sum_i
X_i\otimes Y_i$, then $(\ad X)(x)=\sum_i\hat S (X_i)xY_i$. Thus
\begin{equation} \label{e4.9}
(\ad k)(x)=k^{-1}xk,\ \ (\ad e)(x)=-q^{-1}exk^{-1}+k^{-1}xe,\ \
(\ad f)(x)=-qfxk^{-1}+k^{-1}xf.
\end{equation}

Once we know that $\tilde X^s_j$, $j=-1,0,1$, satisfy
(\ref{e4.10}), then $\tilde X^s_j=\lambda_s X^s_j$ for some
$\lambda_s\in\7C$. This follows immediately from irreducibility of
$U^1$, the fact that $B(H_s)=\oplus^{2s}_{n=0}B(H_s)(n)$ and that
each spectral subspace $B(H_s)(n)$ has multiplicity one. The
constant $\lambda_s$ is uniquely determined by (\ref{e4.7}) or
(\ref{e4.8}). So a direct computation yields $\lambda_s$ as in the
formulation of Lemma.
\ep

Define $\lambda\in M(\hat\A)$ by requiring
$\pi_s(\lambda)=\lambda_s I_s$. It is a straightforward
computation to verify that
$$
\lambda
=\frac{q^2-1}{\sqrt{[2]_q}}(C+\frac{2}{(q-q^{-1})^2}),\ \ \
\hbox{where}\ \
C=fe+\left(\frac{q^\2 k-q^{-\2}k^{-1}}{q-q^{-1}}\right)^2
$$
is the Casimir element. So $\lambda$ belongs to the center of
$\sl2$.

\smallskip

The sphere $(C(S^2_{q,0}),\alpha)$ is a distinguished one. If
$X_{-1}$, $X_0$, $X_1$ are its canonical generators, then the map
$X_j\mapsto -u^1_{j0}$ extends to an embedding of $C(S^2_{q,0})$
into $A$, which intertwines the coaction on~$C(S^2_{q,0})$ with
the left coaction $\D$ of $(A,\D)$ on $A$. Under this embedding
$C(S^2_{q,0})$ is identified with the subalgebra
$$
B=\{a\in A\ |\ f_{it}* a =a\ \forall t\in\7R\}
$$
of $A$ introduced at the beginning of this section. In particular,
the algebra $C(S^2_{q,0})$ also carries a right coaction of $(\hat
A,\hat\D)$. We shall give one more description of $C(S^2_{q,0})$,
which will be pertinent in the sequel.

Let $\Psi$ be the unital C$^*$-subalgebra of $M(\hat A)$ generated
by $\hat A$ and the elements $\lambda^{-1}\chi_j$, $j=-1,0,1$.

\begin{prop}\label{4.5}
\mbox{\ }
\enu{i}
The algebra $\Psi$ has the properties $\Phi (\Psi)\subset
A\otimes\Psi$ and $\hat\D (\Psi)\subset M(\Psi\otimes\hat A)$, so
we have well-defined coactions of $(A,\D)$ and $(\hat A,\hat\D)$
on $\Psi$ and the quotient C$^*$-algebra $\Psi /\hat A$.
\enu{ii}
There exists a unique isomorphism $\sigma\colon
C(S^2_{q,0})\iso\Psi /\hat A$, which maps $X_j$ to
$\lambda^{-1}\chi_j$ ($\mod\hat A$) and intertwines the left
coactions of $(A,\D)$ and the right coactions of $(\hat
A,\hat\D)$.
\end{prop}

\bp
Since $\Phi (\lambda^{-1}\chi_j)=\sum_{k=-1}^1 u^1_{jk}\otimes
(\lambda^{-1}\chi_k)$ by Lemma \ref{4.4}, we obviously have  $\Phi
(\Psi)\subset A\otimes\Psi$. Furthermore, the homomorphism $\Phi$
is a coaction of $(A,\D)$ on $\Psi$ by Corollary \ref{1.4}(i). As
$\pi_s (\lambda^{-1}\chi_j)$, $j=-1,0,1$, are the canonical
generators of $C(S^2_{q,c(s)})$, where
$c(s)=-(q^{2s+1}+q^{-2s-1})^{-2}$, and because $c(s)\to 0$ as
$s\to\infty$, the elements $\lambda^{-1}\chi_j\ \mod\hat A$,
$j=-1,0,1$, satisfy the same relations as the generators of
$C(S^2_{q,0})$. Thus there exists a surjective $*$-homomorphism
$\sigma\colon C(S^2_{q,0})\to \Psi/\hat A$, which maps $X_j$ to
$\lambda^{-1}\chi_j$ (modulo $\hat A$) and intertwines the left
coactions of $(A,\D)$. Since $\Psi/\hat A$ is non-zero, by
Podle\'{s}' classification, this homomorphism must be an
isomorphism.

The assertions that $\hat\D (\Psi)\subset M(\Psi\otimes\hat A)$
and that the isomorphism $\sigma$ intertwines the right coactions
of $(\hat A,\hat\D)$, will be proved in Subsection
\ref{s4.3} below.
\ep

The elements $\chi_j$, $j=-1,0,1$, and $\lambda-\lambda_01$ span
the $4$-dimensional quantum Lie algebra associated to the
$4D_+$-bicovariant calculus of Woronowicz \cite{Wor2,KT}, so
$da=\sum_j(\chi_j*a)\omega_j+((\lambda-\lambda_01)*a)\omega$ in
that context. Thus we can think of $\chi_j$, $j=-1,0,1$, and
$\lambda-\lambda_01$ as left-invariant first order differential
operators, and $\Psi$ as the C$^*$-algebra of left-invariant
pseudodifferential operators of order zero on $SU_q(2)$. Then the
composition $\Psi\to\Psi/\hat
A\stackrel{\sigma^{-1}}{\to}C(S^2_{q,0})$ should be thought of as
the principal symbol.

\smallskip

Before we embark on proving Theorem \ref{4.1}, we end this
subsection by giving some heuristic reasons for why this result
should be true.

The restriction of $P_\phi$ to the center, which is isomorphic to
$c_0 (\7Z_+)$, is not given by a convolution operator on $\7Z$.
However, it is not far from being such an operator. Thus the
theory of random walks on $\7Z$ suggests that the Martin
compactification of the center is obtained by adding one point at
infinity. Now let $H$ be a copy of $H_1$ in $\hat\A$, that is,
there exists a basis $\tilde X_j$, $j=-1,0,1$, in $H$ such that
$\Phi (\tilde X_j)=\sum_{k=-1}^1 u^1_{jk}\otimes \tilde X_k$.
Since $K_{\check\phi}$ commutes with $\Phi$, the elements
$K_{\check\phi}(\tilde X_j)$, $j=-1,0,1$, have the same property
as the elements $\tilde X_j$. Hence, for each $s\in\2\7N$, there
exists a constant $c_s$ such that $K_{\check\phi}(\tilde
X_j)I_s=c_s\lambda^{-1}\chi_j I_s$. The function $\2\7N\ni
s\mapsto c_s$ (which is easily seen to be bounded) embodies
certain properties of the random walk on the center, so it is
natural to assume that it extends to a continuous function on the
Martin compactification. Then the fact that the boundary
consists of one point means that $c_s\to c\in\7C$ as $s\to\infty$.
So, modulo $\hat A$, we have $K_{\check\phi}(\tilde
X_j)=c\lambda^{-1}\chi_j$, and therefore
$K_{\check\phi}(H)\subset\Psi$. Moreover, if $c\ne 0$, then
$\Psi\subset\tilde A_\phi$. This argument can be repeated for all
spectral subspaces (recall that both $C(S^2_{q,0})$ and $\hat A$
have non-zero spectral subspaces for each integer spin), so we get
$K_{\check\phi}(\hat\A)\subset\Psi$. Thus $\tilde
A_\phi\subset\Psi$, and finally $\tilde A_\phi =\Psi$ if at least
one of the functions $s\mapsto c_s$ has a non-zero limit at
infinity.

\subsection{Random Walk on the Center}

The aim of this subsection is to describe the asymptotic behavior
of the Martin kernel on the center. But let us first make some
remarks.

Since $q<1$, we have $n<[n]_q$ for any $n>1$. Hence by
Theorem \ref{2.6}, any state $\phi\in\C$ with $\supp\phi\ne\{0\}$ is
transient. It is worthwhile to note that this result is also valid for
ordinary $SU(2)$. This can be deduced from the existence of
potential kernels for {\it recurrent} random walks on~$\7Z$ (see
e.g. \cite{S}). Nothing like the estimate in Theorem \ref{2.6} is,
however, available. For example, if we consider the random walk
corresponding to the fundamental corepresentation (so $q=1$ and
$\phi=\phi_\2$), then the probability of return to $0$ at the
$n$th step is given by the semicircular law, that is
$p_{\phi^n}(0,0)=\frac{2}{\pi}\int^1_{-1}t^n\sqrt{1-t^2}dt$.

It is known that $U^s\times U^t\simeq \sum^{s+t}_{r=|s-t|}U^r$. It
follows that $\phi$ is generating if and only if
$(\supp\phi)\cap(\2+\7Z_+)\ne\emptyset$.

Since any unitary corepresentation is self-conjugate, we also have
$\phi=\check\phi$ for any positive $\phi\in\C$.

\smallskip

Let now $\phi\in\C$ be a generating positive functional with
$\|\phi\|\le1$. Set
$g_\phi(s,t)=\sum^\infty_{n=0}p_{\phi^n}(s,t)$, so
$G_\phi(I_t)I_s=g_\phi(s,t)I_s$. We want to describe the behavior of
the function $g_\phi(s,t)$ as $s\to\infty$. For this we could
apply the results of Biane \cite{B4} for ordinary $SU(2)$. Namely,
let us for a moment write the superscript '{\it cl}' for the
states on the dual of $SU(2)$. If $\phi=\sum_s\lambda_s\phi_s$, we
set $\phi^{cl}=\sum_s\lambda_s\frac{\dim H_s}{d_s}\phi^{cl}_s$.
Then using Lemma \ref{2.4} and the facts that the fusion
coefficients $N^t_{rs}$ are independent of the deformation
parameter $q\in(0,1]$, it is easy to see that
$p_\phi(s,t)=\frac{\dim H_s}{d_s}p_{\phi^{cl}}(s,t)\frac{d_t}{\dim
H_t}$. In other words, the element $\sum_s\frac{\dim H_s}{d_s}I_s$
is an eigenvector for $P_\phi$ with eigenvalue
$\sum_s\lambda_s\frac{\dim H_s}{d_s}$, and if we consider the Doob
transformation of $P_\phi|_{Z(\hat\A)}$ with respect to this
eigenfunction, we get an operator corresponding to the deformation
parameter $q=1$. Then we can apply the results of Biane to
$\phi^{cl}$ (note that even if $\phi(1)=1$, we have
$\phi^{cl}(1)<1$). We will instead give a slightly different but
more direct proof.

Consider the von Neumann algebra $L(\7T)\cong l^\infty(\7Z)$ of
the circle group $\7T$. Since $\rho$ is an operator
with pure point spectrum $q^\7Z$,
we have an embedding of $L(\7T)$ into $M(\hat A)$ given by
$e^{it}\mapsto u_t=\rho^{-\frac{it}{\log q}}$. Since
$\hat\D(\rho)=\rho\otimes\rho$, this is an embedding of Hopf
algebras. It follows that $P=P_\phi|_{l^\infty(\7Z)}$ is the
operator of convolution with the measure $\phi|_{l^\infty(\7Z)}$.
Let $\{e_n\}_{n\in\7Z}$ be the canonical projections
in~$l^\infty(\7Z)$, so
$$
e_n=\frac{1}{2\pi}\int^{2\pi}_0e^{-int}u_tdt
=\sum_{s=\frac{|n|}{2}}^\infty m^s_{-\frac{n}{2},-\frac{n}{2}}.
$$
If we set $p(n)=\phi(e_n)$, then $P(e_m)e_n=p(m-n)e_n$.

\begin{lem} \label{4.6}
We have $\displaystyle
p_\phi(s,0)=\frac{1}{d_s}(q^{2s}p(-2s)-q^{2s+2}p(-2s-2))$.
\end{lem}

\bp
It is enough to check this for $\phi=\phi_t$. In this case
$P_\phi(I_0)=d_t^{-2}I_t$. So $p_\phi(t,0)=d_t^{-2}$ and
$p_\phi(s,0)=0$ for $s\ne t$. On the other hand,
$p(n)=d_t^{-1}q^n$ if $n\in\{-2t,-2t+2,\dots,2t\}$, and $p(n)=0$
otherwise.
\ep

There exists a unique $\delta_\phi\ge0$ such that
$\phi(\rho^{-\delta_\phi})=1$ (note that the function
$f(t)=\phi(\rho^{-t})$ is convex and $f(-2-t)=f(t)$ as
$\phi_s(\rho^{-t})=d_s^{-1}\Tr\,\pi_s(\rho^{-1-t})$, and
$f(0)=\phi(1)\le1$). Since
$\rho^{-\delta_\phi}e_n=q^{n\delta_\phi}e_n$, this is the same as
to require $\sum_{n\in\7Z}q^{n\delta_\phi}p(n)=1$. Suppose
$\phi=\sum_s\lambda_s\phi_s$. From this point onwards we assume
that
\begin{equation} \label{e4.11}
\sum_{s\in\2\7Z_+}sq^{-2s\delta_\phi}\lambda_s<\infty.
\end{equation}
Equivalently, $\sum_{n\in\7Z}|n|q^{n\delta_\phi}p(n)<\infty$. Then
we set
$$
\lambda_\phi=\sum_{s\in\2\7Z_+}\frac{\lambda_s}{d_s}
   \sum^s_{j=-s}2jq^{2j(1+\delta_\phi)}
=\sum_{n\in\7Z}nq^{n\delta_\phi}p(n).
$$
Note that $\lambda_\phi<0$ as $q<1$.

\begin{prop} \label{4.7}
With notation as above we have
$$
g_\phi(s,0)\sim-\lambda_\phi^{-1}
(1-q^{2+2\delta_\phi})\frac{q^{2s(1+\delta_\phi)}}{d_s}
\ \ \ \hbox{as}\ s\to\infty.
$$
In particular, $\displaystyle
\frac{g_\phi(s+\2,0)}{g_\phi(s,0)}\to q^{2+\delta_\phi}$
as $s\to\infty$.
\end{prop}

The precise values of the constants are unimportant for the
computation of the Martin boundary. What matters is that
$\frac{g_\phi(s+\2,0)}{g_\phi(s,0)}$ converges to some number in
the interval $(0,1)$, and that for states this number is
independent of the particular choice of $\phi$.

\smallskip

\bpp{Proposition \ref{4.7}}
Set $\tilde p(n)=q^{n\delta_\phi}p(n)$ and $\tilde
g(n)=\sum^\infty_{k=0}(\tilde p^{*k})(n)$. Then by the renewal
theorem for $\7Z$ (see \cite[Proposition 24.6]{S}), the sequence
$\{\tilde g(n)\}_n$ tends to $-\lambda_\phi^{-1}$ as
$n\to-\infty$. On the other hand, if
$g(n)=\sum^\infty_{k=0}(p^{*k})(n)$, then $\tilde
g(n)=q^{n\delta_\phi}g(n)$. As
$g_\phi(s,0)=\frac{q^{2s}}{d_s}(g(-2s)-q^2g(-2s-2))$ by
Lemma~\ref{4.6}, we get the desired result (note also that
$d_sd_{s+\2}^{-1}\to q$ as $s\to\infty$).
\ep

\begin{cor}
The Martin compactification of $\2\7Z_+$ with respect to
$P_\phi|_{Z(\hat\A)}$ is obtained by adding one point at infinity.
In particular, if $\phi$ is a state, then the constants are the
only central harmonic elements with respect to $P_\phi$.
\end{cor}

\bp
Since $U^t\times U^r\simeq \sum^{r+t}_{s=|r-t|}U^s$, by Lemma
\ref{2.4}(iv) we get
$P_{\phi_t}(I_r)=\sum^{r+t}_{s=|r-t|}\frac{d_r}{d_td_s}I_s$. As
the algebra $\C$ is commutative in our case, we get
$$
G_\phi(I_t)=d_t^2G_\phi P_{\phi_t}(I_0)=d_t^2P_{\phi_t}G_\phi(I_0)
=d_t^2\sum_{r\in\2\7Z_+}g_\phi(r,0)P_{\phi_t}(I_r),
$$
whence $\displaystyle g_\phi(s,t)=\sum^{s+t}_{r=|s-t|}g_\phi(r,0)
\frac{d_td_r}{d_s}$. By Proposition \ref{4.7} we conclude that for
any $t\in\2\7Z_+$, there exists a constant $c_t$ (depending on
$\phi$) such that $g_\phi(s,t)\sim c_tq^{2s(2+\delta_\phi)}$ as
$s\to\infty$. Thus $\displaystyle
K_\phi(I_t)=\sum_s\frac{g_\phi(s,t)}{g_\phi(s,0)}I_s
=\frac{c_t}{c_0}1\ \mod Z(\hat A)$.
\ep

\subsection{Martin Boundary}\label{s4.3}

Here we will prove Theorem \ref{4.1} and complete the proof of
Proposition \ref{4.5}.

First, analogously to \cite{B4}, we compute the action of the
Martin kernel on certain elements of $\hat\A$.

\begin{prop} \label{4.9}
Let $\phi\in\C$ be a generating positive functional which
has norm not greater than one and satisfies condition (\ref{e4.11}).
Set $\tilde\lambda=\sum_s q^{-2s}I_s$. Then
$K_\phi\F((\alpha^*)^n) =p_n(\tilde\lambda^{-1}k^2)\ \mod\hat A$
for $n\ge0$, where $p_n$ is the polynomial of degree $n$ defined
by the recurrence relation
$$
p_{n+1}(x)=c_\phi p_n(x)x-c_\phi^{-1}p_n(q^{-2}x)(x-1),\ \ p_0=1,
$$
where $c_\phi=q^{2+\delta_\phi}$.
\end{prop}

\bp
We shall prove by induction on $n$ that there exist constants
$a_n(s,m)$, $n\ge0$, $0\le m\le n$, $s\in\2\7Z_+$, such that
$$
G_\phi\F((\alpha^*)^n)=\sum_s\sum^n_{m=0}a_n(s,m)q^{2ms}k^{2m}I_s,
$$
and obtain a recurrence relation for $a_n(s,m)$.

To do this, for $a\in\A$, consider the operator $Q_a$ on $M(\hat A)$
of left convolution with $a$, i.e. $Q_a=(\iota\otimes a)\hat\D$.
Then $P_\phi Q_a=Q_a P_\phi$. As $\F(ab)=Q_a\F(b)$, we have
$P_\phi\F(ab)=Q_aP_\phi\F(b)$ for any
$a,b\in\A$, so $G_\phi\F(ab)=Q_aG_\phi\F(b)$. Hence
\begin{equation} \label{e4.12}
G_\phi\F((\alpha^*)^{n+1})=Q_{\alpha^*}G_\phi\F((\alpha^*)^n)
=\sum_s\sum^n_{m=0}a_n(s,m)q^{2ms}Q_{\alpha^*}(k^{2m}I_s).
\end{equation}
By Lemma \ref{4.3}(ii) and formula (\ref{e4.3}), we have
\begin{equation} \label{e4.13}
k^{2m}I_s=\sum^s_{j=-s}q^{-2mj}m^s_{jj}
 =d_s\sum^s_{j=-s}q^{-2j-2mj}\F(u^s_{-j,-j}).
\end{equation}
Finally, by Lemma \ref{4.3}(iii),
\begin{eqnarray}
Q_{\alpha^*}\F(u^s_{-j,-j})
&=&\F(\alpha^*u^s_{-j,-j})\nonumber\\
&=&d_s^{-1}(q^{-s-j}[s-j+1]_q
\F(u^{s+\2}_{-j+\2,-j+\2})\,
 +\,q^{s-j+1}[s+j]_q\F(u^{s-\2}_{-j+\2,-j+\2}))\nonumber\\
&=&\frac{q^{-s+j-1}}{d_sd_{s+\2}}[s-j+1]_q m^{s+\2}_{j-\2,j-\2}
\,+\,\frac{q^{s+j}}{d_sd_{s-\2}}[s+j]_q m^{s-\2}_{j-\2,j-\2}.
  \label{e4.14}
\end{eqnarray}
Using (\ref{e4.12})--(\ref{e4.14}) a tedious computation yields
\begin{eqnarray*}
a_{n+1}(s,m)
&=&\frac{q^{-2s-1}}{q^{-2s-1}-q^{2s+1}}
      \biggl(-q^{-2(m-1)}a_n(s-\2,m-1)\\
& &\ \ +q^{-2m}a_n(s-\2,m)+a_n(s+\2,m-1)-q^{4s+2}a_n(s+\2,m)\biggr).
\end{eqnarray*}
We also have $a_0(s,0)=g_\phi(s,0)$. Since
$\displaystyle\frac{g_\phi(s+\2,0)}{g_\phi(s,0)}\to
c_\phi=q^{2+\delta_\phi}$ as $s\to\infty$ by Proposition
\ref{4.7}, we see that for any $n$ and $m$, $0\le m\le
n$, there exists a finite limit $\displaystyle
a_n(m)=\lim_{s\to\infty}\frac{a_n(s,m)}{g_\phi(s,0)}$, and these
limits satisfy the recurrence relation
\begin{equation} \label{e4.15}
a_{n+1}(m)=-c_\phi^{-1}q^{-2(m-1)}a_n(m-1)+c_\phi^{-1}q^{-2m}a_n(m)
+c_\phi a_n(m-1).
\end{equation}
Thus if $p_n$ is the polynomial with the coefficients $a_n(m)$, so
$p_n(x)=\sum^n_{m=0}a_n(m)x^m$, then as $\tilde\lambda^{-1}k^2\in
M(\hat A)$ and
$$
K_\phi\F((\alpha^*)^n)=\sum_s\sum^n_{m=0}\frac{a_n(s,m)}{g_\phi(s,0)}
(\tilde\lambda^{-1}k^2)^mI_s,
$$
we conclude that
$K_\phi\F((\alpha^*)^n)=p_n(\tilde\lambda^{-1}k^2)\ \mod\hat A$.
Condition (\ref{e4.15}) can be written as $p_{n+1}(x)=c_\phi
p_n(x)x-c_\phi^{-1}p_n(q^{-2}x)(x-1)$. The leading coefficient of
$p_n$ equals
$$
(c_\phi-c_\phi^{-1}q^{-2(n-1)})(c_\phi-c_\phi^{-1}q^{-2(n-2)})\dots
(c_\phi-c_\phi^{-1})=(-1)^nq^{-n(n-1)}c_\phi^{-n}(c_\phi^2;q^2)_n.
$$
As $c_\phi<1$, it is indeed non-zero, so $p_n$ is a polynomial of
degree $n$.
\ep

\begin{thm} \label{4.10}
Let $\phi\in\C$ be a generating positive functional
which has norm not greater than one and satisfies condition
(\ref{e4.11}). Then
\enu{i}
the Martin compactification $\tilde A_\phi$ coincides with $\Psi$
and, moreover, the linear subspace $K_\phi(\hat\A)+\hat\A$ is
dense in $\Psi$;
\enu{ii}
if $\phi$ is in addition a state, then the Martin kernel $K_\phi$,
regarded as a map $\hat\A\to\Psi/\hat A$, is independent of the
particular choice of $\phi$.
\end{thm}

\bp
By Lemma \ref{4.3} we have $u^s_{ss}=(\alpha^*)^{2s}$ and
$\F(u^s_{ss})=\frac{q^{-2s}}{d_s}m^s_{-s,-s}$. Thus Proposition
\ref{4.9} together with the fact that
$\lambda\tilde\lambda^{-1}\in\7C1+\hat A$, imply
$$
\span\{K_\phi(m^s_{-s,-s})\ |\ s\in\2\7Z_+\}+\hat A
=\7C[\lambda^{-1}k^2]+\hat A,
$$
where $\7C[\lambda^{-1}k^2]$ is the unital algebra generated by
$\lambda^{-1}k^2$. On the other hand, it is easily checked that
\begin{equation} \label{e4.16}
\chi_0=-\lambda+\frac{q\sqrt{[2]_q}}{q-q^{-1}}k^2.
\end{equation}
Thus, we get  $\span\{K_\phi(m^s_{-s,-s})\ |\ s\in\2\7Z_+\}+\hat A
=\7C[\lambda^{-1}\chi_0]+\hat A\subset\Psi$.

Note that the minimal $A$-invariant subspace of $\hat A$
containing $m^s_{-s,-s}$, that is, the subspace spanned by the
elements $(\omega\otimes\iota)\Phi(m^s_{-s,-s})$, $\omega\in A^*$,
coincides with $B(H_s)$. Indeed, this amounts to saying that
$m^s_{-s,-s}$ is a cyclic vector for the adjoint action of $\sl2$
on $B(H_s)$, which in turn can easily be verified using
(\ref{e4.3})--(\ref{e4.5}) and the formulas (\ref{e4.9}) for the
adjoint action. As $\Phi K_\phi=(\iota\otimes K_\phi)\Phi$ by
(\ref{e3.2}), we conclude that $K_\phi(\hat\A)\subset\Psi$.
Furthermore, the linear space $K_\phi(\hat\A)+\hat A$ coincides
with the minimal $A$-invariant subspace of $\Psi$ containing both
$\7C[\lambda^{-1}\chi_0]$ and $\hat A$.

Similarly, if $X_{-1},X_0,X_1$ are the canonical generators of
$B=C(S^2_{q,0})$, then the minimal $A$-invariant subspace of $B$
containing $X^n_0$, for all $n\in\7Z_+$, coincides with the
$*$-subalgebra $\B$ of $B$ generated by $X_j$, $j=-1,0,1$. Indeed,
Podle{\'s} proved \cite{P} that $\B=\oplus^\infty_{n=0}B(n)$, each
spectral subspace $B(n)$ has multiplicity one, and
$(\7C1+B(1))^n\subset\tilde B(n)=\oplus^n_{m=0}B(m)$. Thus to
prove the claim it is enough to check that $X_0^n\in\tilde
B(n)\backslash\tilde B(n-1)$ for any $n\in\7N$. To this end
consider the automorphism $T=(k\otimes\iota)\alpha$ of $\B$. As
$\xi^s_0$ is the only $\pi_s(k)$-invariant vector in $H_s$ by
formula (\ref{e4.3}), the space of $T$-invariant vectors in
$\tilde B(m)$ is $(m+1)$-dimensional. Since $T(X_0)=X_0$, we
conclude that if $X^n_0\in\tilde B(n-1)$, then the elements
$1,X_0,\dots,X^n_0$ are linearly dependent. This is a
contradiction (because e.g. $X_0=\lambda^{-1}\chi_0\ \mod\hat A$).

It follows that modulo $\hat A$ the space $K_\phi(\hat\A)$
coincides with the $*$-algebra generated by $\lambda^{-1}\chi_j$,
$j=-1,0,1$. This proves part (i) of Theorem. To show (ii), note
that by Proposition \ref{4.9} the element $K_\phi(m^s_{-s,-s})\
\mod\hat A$ is independent of the choice of the state $\phi$. As
we explained above, any element in $\hat \A$ is a linear
combination of elements $(\omega\otimes\iota)\Phi(m^s_{-s,-s})$.
Since $K_\phi(\omega\otimes\iota)\Phi(m^s_{-s,-s})
=(\omega\otimes\iota)\Phi K_\phi(m^s_{-s,-s})$, we see that
$K_\phi$, regarded as a map from $\hat\A$ to $\Psi/\hat A$, is
independent of the state~$\phi$.
\ep

Theorem \ref{4.1} would now obviously follow from Theorem
\ref{4.10} and Proposition \ref{4.5}. However, the proof of
Proposition \ref{4.5} is not yet complete. We shall now focus on
how Theorem \ref{4.10} can be used to complete the proof of this
proposition. Due to the fact that $\hat\D(\tilde A_\phi)\subset
M(\tilde A_\phi\otimes\hat A)$, we conclude that
$\hat\D(\Psi)\subset M(\Psi\otimes\hat A)$. Thus there are
well-defined right coactions of $(\hat A,\hat\D)$ on $\Psi$ and
$\Psi/\hat A$. It remains to show that the isomorphism
$\sigma\colon C(S^2_{q,0})\iso\Psi/\hat A$ intertwines the
coactions of $(\hat A,\hat\D)$.

Let $K$ be the map from $\hat\A$ to $\Psi/\hat A$ defined by
$K(x)=K_\phi(x)\ \mod\hat A$, where $\phi$ is any generating state
in $\C$ satisfying condition (\ref{e4.11}). Let $\nu$ be the state
on $\Psi/\hat A=A_\phi$ representing the unit, that is $\nu
K=\hat\psi$ on $\hat\A$. By Proposition \ref{3.12}(ii) and Theorem
\ref{4.10}(i) such a state is unique and $\Phi$-invariant. Since
$C(S^2_{q,0})\subset A$, the state $\varphi|_{C(S^2_{q,0})}$ is the
only invariant state on $C(S^2_{q,0})$. Hence
$\nu\sigma=\varphi|_{C(S^2_{q,0})}$. By Proposition
\ref{1.12}(ii), the state $\varphi$ on $A$ is quasi-invariant
with respect to the right coaction of $(\hat A,\hat\D)$ with
Radon-Nikodym cocycle $y=W(1\otimes\rho^{-2})W^*$. Note that $y\in
M(C(S^2_{q,0})\odot\hat\A)$ as
$$
y=\sum_s\sum_{i,j,l}f_{-2}
   (u^s_{ll})u^s_{il}(u^s_{jl})^*\otimes m^s_{ij}
$$
and $f_z*(u^s_{il}(u^s_{jl})^*)=u^s_{il}(u^s_{jl})^*$. Thus
$\varphi|_{C(S^2_{q,0})}$ is quasi-invariant with Radon-Nikodym
cocycle~$y$. The results of Subsection \ref{s3.3} suggest that
$(K\otimes\iota)\hat\D(I_0)\in M(\Psi/\hat A\odot\hat\A)$ is the
Radon-Nikodym cocycle for $\nu$. Thus if $\sigma\colon
C(S^2_{q,0})\to\Psi/\hat A$ intertwines the right coactions of
$(\hat A,\hat\D)$, we should expect
$(\sigma\otimes\iota)(y)=(K\otimes\iota)\hat\D(I_0)$.

\begin{prop} \label{4.11}
We have $(\sigma\otimes\iota)(y)=(K\otimes\iota)\hat\D(I_0)$.
\end{prop}

\bp
Set $y^s_{ij}=\sum^s_{l=-s}f_{-2}(u^s_{ll})u^s_{il}(u^s_{jl})^*$.
By Lemma \ref{1.2}
$$
(K\otimes\iota)\hat\D(I_0)
 =\sum_s\sum^s_{i,j=-s}K\F(u^s_{ij})\otimes m^s_{ij}.
$$
Thus we must prove that $\sigma(y^s_{ij})=K\F(u^s_{ij})$. We shall
first verify this identity for $i=j=s$. As already remarked in the
proof of Theorem \ref{4.10}, we have $u^s_{ss}=(\alpha^*)^{2s}$.
By Proposition \ref{4.9} it is enough to check that
$y^s_{ss}=p_{2s}(a)$, where $a\in C(S^2_{q,0})$ is determined by
$\sigma(a)=\tilde\lambda^{-1}k^2\ \mod\hat A$, and the polynomials
$p_n$, $n\ge0$, are defined according to the recurrence relation
$$
p_{n+1}(x)=q^2p_n(x)x-q^{-2}p_n(q^{-2}x)(x-1),\ \ p_0=1.
$$
To find the element $a$, note that by definition of $\lambda$
(Lemma \ref{4.4}) and $\tilde\lambda$, we have
$$
\tilde\lambda^{-1}=\frac{1}{(q-q^{-1})\sqrt{[2]_q}}\lambda^{-1}
\ \mod\hat A,
$$
so (\ref{e4.16}) and the definition of $\sigma\colon
C(S^2_{q,0})\to\Psi/\hat A$ yield
$$
\sigma(X_0)=\lambda^{-1}\chi_0\,(\mod\hat A)
=-1+\frac{q\sqrt{[2]_q}}{q-q^{-1}}\lambda^{-1}k^2
 \,(\mod \hat A)=-1+(1+q^2)\tilde\lambda^{-1}k^2\,(\mod\hat A).
$$
Since $X_0=-u^1_{00}$, the element $X_0$ can be found by applying
Lemma \ref{4.3}(iii) to $s=\2$, $j=-\2$. As
$u^\2_{-\2,-\2}=\alpha$, $u^0_{00}=1$ and
$\alpha^*\alpha=1-\gamma^*\gamma$, we get
$X_0=-1+(1+q^2)\gamma^*\gamma$ (though we don't need them, the two
other generators are given by $X_{-1}=\sqrt{q[2]_q}\alpha\gamma^*$
and $X_1=-\sqrt{q[2]_q}\alpha^*\gamma$). Hence $a=\gamma^*\gamma$.

By Lemma \ref{4.3}(iv) we have
$$
y^s_{ss}=\sum^s_{j=-s}f_{-2}(u^s_{jj})u^s_{sj}(u^s_{sj})^*=
\sum^s_{j=-s}q^{-4j}\left[{2s\atop s+j}\right]_{q^2}
(\alpha^*)^{s+j}(\gamma^*\gamma)^{s-j}\alpha^{s+j}.
$$
To prove that $y^s_{ss}=p_{2s}(\gamma^*\gamma)$, we introduce the
polynomials $f_{n,m}$, $n\ge0$, $m\ge0$, by
$$
f_{n,m}(x)=q^{-2mn}x^n(1-x)(1-q^{-2}x)\dots(1-q^{-2(m-1)}x)
 =q^{-2mn}x^n(x;q^{-2})_m.
$$
We assert that
$(\alpha^*)^m(\gamma^*\gamma)^n\alpha^m=f_{n,m}(\gamma^*\gamma)$.
As $\alpha^*\gamma^*\gamma=q^{-2}\gamma^*\gamma\alpha^*$, we need
only to show that $(\alpha^*)^m\alpha^m=(1-\gamma^*\gamma)\dots
(1-q^{-2(m-1)}\gamma^*\gamma)$. Since
$\alpha^*(1-q^{-2l}\gamma^*\gamma)
=(1-q^{-2(l+1)}\gamma^*\gamma)\alpha^*$, this is easily verified
by induction on $m$. It follows that $y^s_{ss}=\tilde
p_{2s}(\gamma^*\gamma)$, where
$$
\tilde p_{2s}(x)=
\sum^s_{j=-s}q^{-4j}\left[{2s\atop s+j}\right]_{q^2}f_{s-j,s+j}(x).
$$
It remains to prove that $\tilde p_n=p_n$, or equivalently, that
the polynomials $\tilde p_n$, $n\ge0$, satisfy the relation
$$
\tilde p_{2s}(x)=q^2\tilde p_{2(s-\2)}(x)x
 -q^{-2}\tilde p_{2(s-\2)}(q^{-2}x)(x-1).
$$
Using the identities
$$
\left[{2s\atop s+j}\right]_{q^2}=q^{2(s+j)}
  \left[{2s-1\atop s+j}\right]_{q^2}
 +\left[{2s-1\atop s+j-1}\right]_{q^2},
$$
and $q^{2(s+j)}f_{s-j,s+j}(x)=f_{s-j-1,s+j}(x)x$ and
$f_{s-j,s+j}(x)=f_{s-j,s+j-1}(q^{-2}x)(1-x)$, we see that the
polynomial $\tilde p_{2s}$ equals
$$
\sum^{s-1}_{j=-s}q^{-4j}\left[{2s-1\atop s+j}\right]_{q^2}
 f_{s-j-1,s+j}(x)x+\sum^s_{j=-s+1}q^{-4j}\left[{2s-1\atop
s+j-1}\right]_{q^2}f_{s-j,s+j-1}(q^{-2}x)(1-x).
$$
The first summand is $\displaystyle\
q^2\hspace{-3mm}\sum^{s-\2}_{j=-s+\2}q^{-4j}\left[{2(s-\2)\atop
s-\2+j}\right]_{q^2}f_{s-\2-j,s-\2+j}(x)x=q^2\tilde
p_{2(s-\2)}(x)x$, whereas the second one is $q^{-2}\tilde
p_{2(s-\2)}(q^{-2}x)(1-x)$. Thus $\tilde p_n=p_n$ and we have
demonstrated the equality $\sigma(y^s_{ss})=K\F(u^s_{ss})$.

\smallskip

Consider now the linear maps $T_1\colon B(H_s)\to C(S^2_{q,0})$
and $T_2\colon B(H_s)\to \Psi/\hat A$ given by
$$
T_1(m^s_{ij})=(-1)^{i-j}q^{-i-j}d_sy^s_{-i,-j}\ \ \hbox{and}\ \
T_2(m^s_{ij})=(-1)^{i-j}q^{-i-j}d_sK\F(u^s_{-i,-j}).
$$
By Lemma \ref{4.3}(ii) the map $T_2$ is just $K|_{B(H_s)}$, so it
intertwines the left coactions of $(A,\D)$, i.e. $\Phi
T_2=(\iota\otimes T_2)\Phi$. By the same lemma we also have
$$
(-1)^{i-j}q^{-i-j}d_sy^s_{-i,-j}=(-1)^{i-j}q^{-i-j}d_s
 \sum_lq^{4l}u^s_{-i,-l}(u^s_{-j,-l})^*
=d_s\sum_lq^{2l}(u^s_{il})^*u^s_{jl}.
$$
Since
$$
\D((u^s_{il})^*u^s_{jl})=\sum_{n,m}(u^s_{in})^*u^s_{jm}
 \otimes(u^s_{nl})^*u^s_{ml}
\ \ \hbox{and}\ \
\Phi(m^s_{ij})={U^s}^*(1\otimes m^s_{ij})U^s
=\sum_{n,m}(u^s_{in})^*u^s_{jm}\otimes m^s_{nm},
$$
we see that $T_1$ also intertwines the left coactions of $(A,\D)$.
Therefore the maps $\sigma T_1$ and $T_2$ intertwine the coactions
of $(A,\D)$ and coincide on the element $m^s_{-s,-s}$. Since the
minimal $A$-invariant subspace of $B(H_s)$ containing
$m^s_{-s,-s}$ is the whole of $B(H_s)$, we get $\sigma T_1=T_2$.
Hence $\sigma(y^s_{ij})=K\F(u^s_{ij})$ for all $i$ and $j$.
\ep

The following corollary completes the proof of Proposition
\ref{4.5}, and thus also of Theorem \ref{4.1}.

\begin{cor} \label{4.12}
The isomorphism $\sigma\colon C(S^2_{q,0})\iso\Psi/\hat A$
intertwines the right coactions of $(\hat A,\hat\D)$.
\end{cor}

\bp
The elements $y=W(1\otimes\rho^{-2})W^*\in
M(C(S^2_{q,0})\odot\hat\A)$ and $\tilde
y=(K\otimes\iota)\hat\D(I_0)\in M(\Psi/\hat A\odot\hat\A)$ are
invertible cocycles for the right coactions of $(\hat A,\hat\D)$.
As $(\sigma\otimes\iota)(y)=\tilde y$ we get
$$
(\sigma\otimes\iota\otimes\iota)(\hat\Phi\otimes\iota)(y)
=(\sigma\otimes\iota\otimes\iota)
  ((\iota\otimes\hat\D)(y)(y^{-1}\otimes1))
=(\iota\otimes\hat\D)(\tilde y)(\tilde y^{-1}\otimes1)
=(\hat\D\otimes\iota)(\tilde y).
$$
Thus $(\sigma\otimes\iota)\hat\Phi(a)=\hat\D\sigma(a)$ for any
element $a$ of the form $(\iota\otimes\omega)(y)$, $\omega\in\A$,
and therefore also for any element in the C$^*$-algebra generated
by $(\iota\otimes\omega)(y)$, $\omega\in\A$. Since
$(\iota\otimes\omega)(\tilde y)=K(\iota\otimes\omega)\hat\D(I_0)$,
the set of elements $(\iota\otimes\omega)(\tilde y)$,
$\omega\in\A$, is identical to $K(\hat\A)$. The set $K(\hat\A)$ is
dense in $\Psi/\hat A$ due to Theorem \ref{4.10}(i). Since
$\sigma$ is an isomorphism, the set of elements
$(\iota\otimes\omega)(y)$, $\omega\in\A$, is dense
in~$C(S^2_{q,0})$ as well. So $\sigma$ intertwines the right
coactions of $(\hat A,\hat\D)$.
\ep

Now since $\nu\sigma=\varphi$, the state $\varphi|_{C(S^2_{q,0})}$
is quasi-invariant with Radon-Nikodym cocycle
$y=W(1\otimes\rho^{-2})W^*$ and
$(\sigma\otimes\iota)(y)=(K\otimes\iota)\hat\D(I_0)$, we obtain
the following result.

\begin{cor} \label{4.13}
The state $\nu$ on $\Psi/\hat A$ is quasi-invariant with respect
to the right coaction of $(\hat A,\hat\D)$ with Radon-Nikodym
cocycle $(K\otimes\iota)\hat\D(I_0)$.
\epp
\end{cor}

\subsection{Poisson Boundary}

Having the computations of the previous subsection behind us, it
is now an easy matter to determine the Poisson boundary. By
Theorem
\ref{4.10} we can identify the Martin boundary with $\Psi/\hat A$.
Moreover, we know already that the state $\nu=\varphi\sigma^{-1}$
on $\Psi/\hat A$ represents the unit, $K=K_\phi\ \mod\hat A$ is
independent of~$\phi$, and $K(\hat\A)$ is dense in $\Psi/\hat A$.
By Proposition \ref{3.12}(ii) the map $K^*\colon\pi_\nu(\Psi/\hat
A)''\to M(\hat A)$ defined by $(K(x),a)_\nu=(x,K^*(a))_{\hat\psi}$
for $x\in\hat\A$ and $a\in\pi_\nu(\Psi/\hat A)''$, is an
isomorphism of $\pi_\nu(\Psi/\hat A)''$ onto $H^\infty(M(\hat
A),P_\phi)$ which intertwines the left coactions of $(M,\D)$.
Since $\nu$ is quasi-invariant with respect to the right coaction
of $(\hat A,\hat\D)$, by the discussion following the proof of
Proposition~\ref{3.11}, the right coaction of $(\hat A,\hat\D)$
extends to a right coaction of $(\hat M,\hat\D)$ on the von
Neumann algebra~$\pi_\nu(\Psi/\hat A)''$. Then by Proposition
\ref{3.11}(ii) and Corollary~\ref{4.13}, the map $K^*$ intertwines
the right coactions of $(\hat M,\hat\D)$, and
$K^*=(\nu\otimes\iota)\hat\D$. We summarize this discussion in the
following theorem which makes the statement of Theorem \ref{4.2}
more precise.

\begin{thm} \label{4.14}
\mbox{\ }
\enu{i}
For any generating state $\phi=\sum_s\lambda_s\phi_s\in\C$ such
that $\sum_s s\lambda_s<\infty$, the map
$K^*\colon\pi_\nu(\Psi/\hat A)''\to M(\hat A)$ gives an
isomorphism of $\pi_\nu(\Psi/\hat A)''$ onto $H^\infty(M(\hat
A),P_\phi)$ intertwining the left coactions of $(M,\D)$ and the
right coactions of $(\hat M,\hat\D)$. In particular,
$H^\infty(M(\hat A),P_\phi)\subset M(\hat A)$ is independent
of~$\phi$.
\enu{ii}
We have $K^*=(\nu\otimes\iota)\hat\D$, so
$K^*\sigma=(\varphi\otimes\iota)\hat\Phi$ on $C(S^2_{q,0})\subset A$.
\epp
\end{thm}

\bigskip

\flushleft
{Sergey Neshveyev\\
Mathematics Institute\\
University of Oslo\\
PB 1053 Blindern\\
Oslo 0316\\
Norway\\
{\it e-mail}: neshveyev@hotmail.com}

\flushleft
{Lars Tuset\\
Faculty of Engineering\\
Oslo University College\\
Cort Adelers st. 30\\
Oslo 0254\\
Norway\\
{\it e-mail}: Lars.Tuset@iu.hio.no}

\end{document}